\documentclass[a4paper, 12pt]{article}
\usepackage[english]{babel}

\usepackage{amsthm, amssymb}
\usepackage[latin1]{inputenc}
\usepackage[T1]{fontenc}

\newtheorem{hypothesis}{Hypothesis}[section]
\newtheorem{proposition}{Proposition}[section]
\newtheorem{theorem}{Theorem}[section]
\newtheorem{corollary}{Corollary}[section]

\newtheorem{remark}{Remark}[section]
\newtheorem{lemma}{Lemma}[section]
\newcommand{\eqref}[1]{(\ref{#1})}
\newcommand{\ptl}{\partial}

\newcommand{\RR}{\mathbb{R}}
\newcommand{\NN}{\mathbb{N}}
\newcommand{\refeq}[1]{(\ref{#1})}
\newcommand{\calX}{{\cal X}}
\newcommand{\calJ}{{\cal J}}
\newcommand{\norm}[1]{\left\vert#1\right\vert}
\newcommand{\Norm}[1]{\left\Vert#1\right\Vert}

\begin{document}

\title{Asymptotic Expansions for the Heat Kernel and the Trace of a Stochastic Geodesic
Flow}

\bibliographystyle{plain}

\author{\small Sergio Albeverio\\
\small Institut f. Angew. Mathematik\\
\small Universit\"at Bonn \\
\small Sergio.Albeverio@uni-bonn.de\\
\small CERFIM, Locarno
\and
\small Astrid Hilbert\\
\small MSI\\
\small Växjö Universitet, Sweden\\
\small Astrid.Hilbert@vxu.se \\
\and
\small Vassily Kolokoltsov \\
\small Department of Statistics\\
\small University of Warwick\\
\small v.Kolokoltsov@warwick.ac.uk}
\maketitle

\begin{abstract}
We analyze the asymptotic behaviour of the heat kernel defined by a stochastically
perturbed geodesic flow on the cotangent bundle of a Riemannian manifold
for small time and small diffusion parameter.
This extends WKB-type methods to a particular case of a degenerate Hamiltonian.
We derive uniform bounds for the solution of the degenerate Hamiltonian
boundary value problem for small time. From this equivalence of solutions of the Hamiltonian equations
and the corresponding Hamilton Jacobi equation follows. The results are exploited to
derive two sided estimates and multiplicative asymptotics for the heat kernel and the
trace.
\end{abstract}

\section{Introduction}

The problem of finding Gaussian bounds for the heat kernel and of deriving asymptotic
expansions for non degenerate diffusions has been studied by various researchers.
Expansions for the stochastic flow and stochastic variational calculus have been exploited by
 ~\cite{Mol}, \cite{Dav88}, \cite{AlBr}, \cite{KuSt87}, \cite{BA89}, and \cite{Lea94}.
For hypoelliptic operators satisfying a uniform weak Hörmander hypothesis uniform in
space estimates for the heat kernel have been given, e.g., by L{\'e}andre. For
Gaussian hypoelliptic generators D. Manankiandrianana~\cite{ChE81} and M. Chaleyat,
L. Elie~\cite{Man79} have classified the structure of the corresponding processes,
derived the associate action functionals, and the associated subriemannian metrics.

In \cite{Mas85} Maslov observed that for a class of pseudodifferential equations
(called tunnel equations) on $\RR^n$ which includes the example of elliptic
second order equations asymptotic bounds for the heat kernel can be established
for small values of a parameter $h$ (Planck's constant). Results using WKB-type methods
may be found e.g. in \cite{Mas72}, \cite{Mas85}, \cite{MaF81}, \cite{KoM97},
\cite{Kol00}.
In this paper we modify the WKB-type method by introducing a change of variable which
adjust the different time scales of the variables. We derive bounds for the kernel
of degenerate diffusions associated with the geodesic flow on a compact $d$-dimensional
Riemannian manifold $M$ with metric g which is subject to a stochastic perturbation in
the tangent space variable.
In local coordinates the diffuion on the cotangent bundle $T^*M$ is given by.

\begin{eqnarray}\label{Greensfunktion}
  \frac {\partial u}{\partial t} &=&  {\cal L}u
  = \frac h2 \left(g_{ij}(x)\frac {\partial^2 u}{\partial y_i \partial y_j}\right)
  + \left(G(x)y, \frac {\partial u}{\partial x}\right)
  - \frac 12\left(\frac {\partial }{\partial x}(G(x)y,y),\frac {\partial u}{\partial y}\right)\\
& = &: h^{-1} H\left(x,y,-h\frac \partial{\partial x}, -h\frac \partial{\partial y}\right)u \nonumber
\end{eqnarray}
with initial data $u(0,x,y,x_0,y^0) = \delta(x-x_0) \delta(y-y^0)$ for $x$ and $y$ in a
neighbourhood of the point $(x_0,y^0)\in T^*M$, $h>0$ is a (small) parameter $g_{ij}(x)$
is the expression of the metric in local coordinates with inverse $G^{ij}(x)$. The
summation convention over indices is used. The pseudodifferential operator
$H(x,y,-\frac {\partial }{h\partial x},-\frac {\partial}{h\partial y})$ defines a
Hamiltonian function $H(x,y,p,q)$
the
arguments of which can be divided into the variables $x,y \in \mathbb R^d$ and their
respective  conjugate variables i.e. the momenta $p$ and $q$.
In local coordinates the Hamiltonian takes the form
\begin{equation}\label{Hamilton}
  H(x,y,q,p) = \frac 12 \left(g(x)q,q\right)-\left(G(x)y,p\right)
       +\frac 12\left(\frac {\partial}{\partial x}(G(x)y,y),q\right)\ .
\end{equation}
It is possible to write a stochastic differential equation associated with the operator
$\cal L$. One way is to introduce an imbedding $r$ of the manifold $M$ into a large
Euclidian space $\mathbb{R}^N$ (such an imbedding always exists).
In local coordinates the corresponding stochastic system will have the form
\begin{eqnarray}\label{eq:5}
  dx_i &=& \frac {\partial H}{\partial y_i}dt \nonumber \\
  dy_i &=& -\frac {\partial H}{\partial x_i}dt
          + \frac {\partial}{\partial x_i}(\gamma, d\omega)
\end{eqnarray}
where $(x.y)\in T^*M$, $\omega$ is standard Brownian Motion in $\mathbb{R}^d$ and
$(\cdot,\cdot)$ denotes the standard scalar product in $\mathbb{R}^N$. The system
describes a Newtonian particle corresponding to a Hamiltonian system perturbed by
a Gaussian White Noise force $``\dot{w}``$ with the friction
$\frac 12 \gamma_i^{kl}y_k y_l$ see e.g. \cite{Nel67},
\cite{AHZ92}, \cite{MaS00} for the flat
case or \cite{Joe78} for a covariant definition. In the sequel neglecting friction
we shall restrict ourselves to the example of the geodesic flow on the cotangent
bundle $T^*M$ of a Riemannian manifold $M$. In this case
the Hamiltonian corresponding to the deterministic part of \eqref{eq:5} is of the
form $f(x) = \frac 12 (G(x)y,y)$ where $G(x) = g^{-1}(x)$ is smooth by definition,
see also \refeq{Hamilton}.

\vskip 0.3cm

The content of the paper is organized in five sections. The second outlines
the general  structure of an invariant diffusion. In the subsequent one we study the
underlying Hamiltonian dynamics of the stochastic geodesic flow. This results in the
uniform estimate \refeq{eq:2.3} and in the existence of a local diffeomorphic flow on
the cotangent bundle. In Section 4 we explain the WKB-method for imaginary time and
the change of variables which compensates for the degeneracy and allows for a closed
system of algebraic equations for the coefficients of an asymptotic expansion for the
solution of \refeq{Greensfunktion} in powers of $h$. The remaining part of this section is dedicated
to the existence of a multiplicative error. The paper is concluded by discussing
an asymptotic expansion for the trace formula for the generator corresponding to a
stochastic geodesic flow \refeq{Greensfunktion}. The coefficients of the expansion are
geometric invariants of the manifold.

\section{Geometric background for the system under consideration}

Let us introduce a diffusion process on the cotangent bundle of a compact Riemannian
manifold $M$ with metric $g$ which in local coordinates has the form

\begin{eqnarray}\label{eq:1}
 h\frac \partial{\partial t} u
  &=& \frac {h^2}2 g_{ij}(x) \frac {\partial^2}{\partial y_i \partial y_j} u
      + h(a^i(x) + \alpha^{ij}(x)y_j)\frac {\partial u}{\partial x^i} \nonumber\\
  & & +~h(b_i(x) + \beta_i^j(x)y_j
      + \frac 12 \gamma_i^{jl}(x)y_jy_l)\frac {\partial u}{\partial y_i} - V(x,y)u
\end{eqnarray}
where $h$ is a strictly positive parameter, the $d\times d$-matrix $g(x) = (g_{ij}(x))$
is positive definite, $\alpha,\ \beta,\ \gamma$ are $d\times d$-matrices, $\alpha$ being
non-degenerate, $a,b$ are uniformly bounded in $x\in\mathbb{R}^d$ and $V$ is a polynomial
of degree 4 in $y \in \mathbb{R}^d$ with coefficients depending on $x\in\mathbb{R}^d$
bounded from below. This diffusion defines a regular degenerate diffusion of rank one
in the sense of \cite{Kol00}. In this book there are given necessary and sufficient
conditions for this diffusion to be invariant.

Let us show that the pseudodifferential operator in \eqref{Greensfunktion} is an
invariant degenerate diffusion. In order to study invariance under change of coordinates
in the manifold, we recall that changing the
variable $x \rightarrow \tilde x$ for the diffusion \eqref{Greensfunktion}
induces a change of the momentum variables according to the rule
$\tilde{y} = y\frac {\partial x}{\partial \tilde{x}}$. This implies e.g.
\begin{eqnarray}\label{eq:3}
  \frac {\partial \tilde{y}_k}{\partial y_j} = \frac {\partial x^j}{\partial \tilde{x}^k},\quad
  \frac {\partial u}{\partial x^i}
   = \frac {\partial u}{\partial \tilde{x}^{\ell}}\frac {\partial \tilde{x}^{\ell}}{\partial x^i}
   + \frac {\partial u}{\partial \tilde{y}_{\ell}}\frac {\partial \tilde{y}_{\ell}}{\partial x^i}\\
  \frac {\partial u}{\partial y_i}
  = \frac {\partial u}{\partial \tilde{y}_{\ell}} \frac {\partial \tilde{y}_{\ell}}{\partial y_i},
  \quad \frac {\partial^2 u}{\partial y_i \partial y_j}
  = \frac {\partial^2 u}{\partial \tilde{y}_k \partial \tilde{y}_{\ell}}
  \frac {\partial \tilde{y}_{\ell}}{\partial y_j}\frac {\partial \tilde{y}_k}{\partial y_i} .
\label{eq:4}
\end{eqnarray}
for $i,j,k=1,\ldots ,d$. An immediate consequence of \eqref{eq:3} and \eqref{eq:4} is
that the terms of different orders do not mix and hence may be considered separately.
Then the invariance of equation \eqref{eq:1} follows from direct calculation
using the relations \eqref{eq:3} and \eqref{eq:4} and the fact that the metric tensor
transforms according to
$\tilde{g}_{ij}(\tilde{x}) = g_{ij}(x)\frac {\partial x^j}{\partial \tilde{x}^{\ell}}
  \frac {\partial x^i}{\partial \tilde{x}^k}$.

The calculations in this paper will be carried out in normal coordinates around
$x_0$. We are using the definition given in \cite{CFKS87}
for which $x=0$ and for which the matrix associated with the
Riemannian metric satisfies $\det g(x) = 1$, which can be achieved by a
change of variables and the expansion
\begin{equation}\label{eq:8}
  g_{ij}(x) = \delta_i^j + \frac 12 g_{ij}^{kl}x^k x^l + O(|x|^3).
\end{equation}
These conditions imply
\begin{displaymath}
  \sum_{i=1}^n g_{ii}^{kl} = 0 \quad\mbox{for}\quad 1\leq k,l \leq d
\end{displaymath}
and the following representation for the Gausssian (scalar) curvature
\begin{displaymath}
  R=\sum_{i,k} g_{ik}^{ik}.
\end{displaymath}

For the corresponding inverse matrix $G(x) = g^{-1}(x)$ and its derivative we get the
following Taylor approximations

\begin{equation}\label{eq:9}
  G^{ij}(x) = \delta_i^j - \frac 12 g^{ij}_{kl} x^k x^l + O(|x|^3)
  \quad\mbox{and}\quad
    \Gamma^{ij}_k(x) = -g_{kl}^{ij}x^{\ell} + O(|x|^3).
\end{equation}

Here $\Gamma^{ij}_k(x):= \frac {\partial G_{ij}}{\partial x^k}(x)$ are the coefficients
of the Levy Civita connection in local coordinates.
For symmetry reasons we also introduce the notation
$\gamma_{i\ell j}=\frac{\partial{g_{i\ell}}}{\partial{x^j}}$ which will be used in
the next section.

\section{Hamiltonian Dynamics}
The Hamiltonian system corresponding to \eqref{Greensfunktion} is given by
\begin{eqnarray}\label{eq:2.1}
 \dot{x}^i & = & \frac{\partial H}{\partial p_i} = -G^{ij}(x) y_j\\
 \dot{y}_i & = & \frac{\partial H}{\partial q^i} =
                 \frac{1}{2}\frac{\partial}{\partial x^i}\left(G(x)y,y\right)
                 + \frac{1}{2} g_{ij}(x)q^j
          = \frac{1}{2} \Gamma_i^{kj}(x) y_k y_j + \frac{1}{2} g_{ij}(x)q^j  \nonumber\\
 \dot{q}^i & = & -\frac{\partial H}{\partial y_i} =
                 G^{ij}(x)p_j - \left(\frac{\partial}{\partial x}G^{ij}(x)y_j,q\right)
          = G^{ij}(x)p_j - \Gamma_k^{ij}(x) y_j q^k
               \nonumber\\
 \dot{p}_i & = & -\frac{\partial H}{\partial x^i} =
               \left(\frac{\partial}{\partial x^i}G(x)y,p\right)
        - \frac{1}{2}\left(\frac{\partial}{\partial x^i}\frac{\partial}{\partial x}\left(G(x) y,y\right),q\right)
        -\frac{1}{2}\left(\frac{\partial}{\partial x^i}g(x) q,q\right)
        \nonumber
\end{eqnarray}
with initial condition $x(0)=x_0,\ y(0)=y^0,\ q(0)=q_0,\ p(0)=p^0$. The
corresponding flow will be denoted by $(X,Y,Q,P)$.
\begin{hypothesis}\label{hypo:H}
Let there exist constants $k, t_0$ and $c_0$ such that for all $c \in (0,c_0]$ and
$t \in (0, t_0]$ the solution to (\ref{eq:2.1}) exists on $[0,t]$ whenever the
initial data $y^0,q_0,p^0$ are satisfying
$|y^0|\leq\frac{c}{t},~ | q_0|\leq\frac{c^2}{t^2},~ |p^0| \leq \frac{c^3}{t^3}$
and $x_0$ may be chosen freely.
\end{hypothesis}
\begin{proposition}\label{prop:2.1}
Under the hypothesis \ref{hypo:H} the components of the flow
$(X,Y,P,Q)$ started at $(x_0,y^0,p^0,q_0)$ satisfy
\[ |x(t)-x_0|\leq k t(1+\frac{c}{t}),\quad |y(t)-y^0|\leq k t(1+\frac{c^2}{t^2})\]
\[ |q(t)-q_0|\leq k t\left(1+\frac{c^3}{t^3}\right),\quad |p(t)-p^0|
      \leq k t(1+\frac{c^4}{t^4}).\]
\end{proposition}
\begin{proof}
The proof uses Taylor expansion and is a direct consequence of the fact that the
matrix $g$ is a bounded function of the variable $x$
(c.f. also Proposition 2.3.1 in \cite{Kol00}).
\end{proof}

If we associate with the variables $x,y,q,p$ the degree $0,1,2,3$, respectively each
right hand side of \refeq{eq:2.1} is of specific degree varying between 1 for
$\dot{x}$ and 4 for $ \dot{p}$. The degree is related to the power of $\frac 1t$
which is crucial for the iteration procedure giving the expansions intrinsic to
the WKB-method.
\vskip 0.2cm

In the subsequent proposition we show that the solution $(x,y,q,p)$ of
\refeq{eq:2.1} starting at $(x_0,y^0,q_0,p^0)$ can be represented
in terms of the solution $(\tilde x, \tilde y, \tilde p, \tilde q )$ solving
\refeq{eq:2.1} with initial value $(x_0,y^0,0,0)$.
For reason of space we shall skip the arguments of the functions $g(x),\ G(x)$ and
of their derivatives in future.

\begin{proposition}\label{prop:2.2}
Under the assumptions of Proposition~\ref{prop:2.1} we have:
\begin{eqnarray*}
  x & = & \tilde x - \frac{1}{2}t^2q_0 - \frac{1}{6}t^3[p^0+ \left( \Omega'y^0 \right) q_0]
                   + \delta^4 \\ \nonumber
  y & = & \tilde y + tq_0 + \frac{1}{2}t^2[p^0 + \left( \Omega y^0 \right) q_0]
                   + \frac{1}{t}\delta^4 \\\nonumber
  q & = & \tilde q + t[G(x_0)p^0+(\Gamma y^0)q_0] \delta^4
                                + \frac 12 t^2[(\Lambda y^0)p^0+(\Lambda'(y^0)^2)q_0
                                +\Lambda''(q_0)^2] + \delta\\\nonumber
  p & = & \tilde p + t \delta^4.
\end{eqnarray*}
where the error term $ \delta $ is of the order $ O(t+c)$. Moreover,
$\tilde x = x_0 - y_0t + O(t^2)$, $\tilde y = y^0 + O(t)$, $\tilde p = O(t)$,
$\tilde q = O(t)$, and
\begin{eqnarray*}
\left(\Gamma y_0 \right)_{i}^k    &=& \Gamma_i^{jk} y_j^0
              = \frac{\partial G^{j k}}{\partial x^i}y_j^0, \qquad\qquad
 \left(\tilde\Gamma y^0,y^0\right)_{ik} = \tilde\Gamma_{ik}^{\ell j} y^0_\ell y^0_j
         = \frac{\partial{\Gamma_{i}^{\ell j}}}{\partial{x^k}}  y^0_\ell y^0_j\\
 \left(\Omega y_0 \right)_{i\ell} &=& \left[ g_{ij} \Gamma^{jm}_{\ell} g_{j\ell}
               -\gamma_{i\ell j} G^{jm} \right](x_0) y^0_m ,\\
 \left(  \Omega' y^0 \right)_{ik} &=& \Gamma_i^{j\ell}g_{jk}y_\ell^0
                - \gamma_{ikm}G^{m\ell}y_\ell^0 - g_{ij}\Gamma_k^{j\ell}y_\ell^0\\
 \left(\Lambda y_0 \right)^{ik}    &=& -\Gamma_m^{ik}G^{m\ell} y_\ell^0
      + G^{im} \Gamma_m^{k\ell}  y_\ell^0 -  G^{mk} \Gamma_m^{ij}  y_j^0 ,\\
 \left(\Lambda' (y_0)^2 \right)_{k}^i  &=&-\frac 12 G^{im}\tilde\Gamma_{mk}^{\ell j} y_\ell^0 y_j^0
      - \tilde\Gamma_{mk}^{ij} G^{m\ell} y_\ell^0 y_j^0
      -\frac 12 \Gamma_k^{ij}\Gamma_j^{m\ell} y_\ell^0 y_m^0
      -\Gamma_m^{ij}\Gamma_k^{m\ell} y_j^0y_\ell^0 .\\
\left(\Lambda''\right)^i_{k\ell}    &=&  - \frac 12  G^{im}\gamma_{k\ell m}
               - \frac 12 g_{m\ell}\Gamma_k^{im}
\end{eqnarray*}
\end{proposition}

\begin{proof}
For the proof we expand the components
$p-\tilde p,\, q-\tilde q,\, y-\tilde y,\, x-\tilde x$
in Taylor series with respect to time. For the time derivatives of the expansion
we insert the explicit expressions given by the system of Hamilton equations
(\ref{eq:2.1}).
Due to linearity of differentiation it suffices to study the derivatives of $y(t)$
and insert $q_0=0$ and $p^0=0$ to retrieve the derivatives of $\tilde y(t)$.
We begin with the function $\Delta p_i(t)= p_i(t)-\tilde p_i(t)$
because it is of highest degree and hence determines the degree to be considered
for all other components. We have
$\Delta p_i(t)=p_i^0 + \int_0^t \dot{\Delta p}_i(\tau) \,d\tau$. We directly find:
\begin{eqnarray} \label{eq:2.1a}
 \Delta p_i(t) &=& p_i^0 + \int_0^t \Gamma_i^{k\ell} y_\ell p_k
            - \frac 12  \tilde\Gamma_{mi}^{k\ell} y_\ell y_m q^k
            - \frac 12 \gamma_{ijk} q^j q^k \, d\tau\quad \\
          &=& p_i^0 + \int_0^t L(y^0p^0,(y^0)^2q_0) + O(\tau^2) \, d\tau \ .
\end{eqnarray}
Here $L$ is a linear function of degree $4$ in its arguments with bounded
coefficients depending on $x(t)\in\RR^d$. For $\dot{\tilde p}$ the term  corresponding
to $L$ vanishes except for the uniform error $O(t^3)$.
Inserting the estimates given in Proposition \ref{prop:2.1} in the integrand
gives the desired result.
For the function
$\Delta q_i(t)= q_i(t)-\tilde q_i(t)= q^i_0 + t \dot{\Delta q^i}_0
    +\int_0^t (t-\tau) \ddot{\Delta q^i}(\tau) \,d\tau$ we set off analogously:
\begin{eqnarray}
\dot{q^i}(\tau)  &=& G^{ij} p_j - \Gamma_\ell^{ij}y_j,q^\ell
                    \ \mbox{in particular}\ \dot{q^i_0}(\tau)
   = G^{ij} p_j^0 - \Gamma_\ell^{ij}y_j^0,q^\ell_0\nonumber\\
\ddot{q^i}(\tau) &=& \frac{\partial G^{ij}}{\partial x} \dot x y_j + G^{ij} \dot p_j
         - \left(\frac{\partial^2 G^{ij}}{\partial x \partial x} \dot x y_j,q\right)
         - \left(\frac{\partial G^{ij}}{\partial x}\dot{y}_j,q\right)
         - \left(\frac{\partial G^{ij}}{\partial x} y_j,\dot{q}\right) \nonumber
\end{eqnarray}
Direct calculation gives a long expression the structure of which becomes clear
by inserting the notation in the proposition.
\begin{displaymath}
 \ddot{q}(\tau)   = (\Lambda y)p + (\Lambda' y^2)q + \Lambda'' q^2
         = (\Lambda y)p^0 + (\Lambda' y^2)q_0 + \Lambda'' q^2_0 + O(t)\ .
\end{displaymath}

A direct consequence is that
$\dot{\tilde q^i(0)}$ and $\ddot{\tilde q^i(0)}$ vanish except for an error of
order $O(t)$.
We find using the estimates of Proposition \ref{prop:2.1} under the integral
\begin{eqnarray}
\Delta q^i(t) &=& q^i_0 + t L(y^0q_0,p^0)
          + \int_0^t (t-\tau) L(y(\tau)p(\tau),y^2(\tau)q(\tau),q^2(\tau)) \,d\tau\nonumber\\
       &=& q^i_0 + t L(y^0q_0,p^0)
                         + t^2 L(y^0p^0,(y^0)^2q_0,q^2_0) + O(t^3)\label{eq:2.1b}
\end{eqnarray}
where we proceeded as above.
In order to achieve an analogous representation for the function $y(t)$ the
starting equation is given by
\begin{equation} \label{eq:2.1c}
 \Delta y_i(t) =y_i(t)-\tilde y_i(t)= y_i^0 + t \dot{\Delta y_i}^0
               + \frac 12 t^2 \ddot{\Delta y_i}^0
  +\int_0^t (t-\tau)^2 \Delta (y_i)^{(3)}(\tau) \,d\tau \ .
\end{equation}
The first order derivatives are explicitly given in \refeq{eq:2.1}.
For $2\ddot y_i(t)$ we have
\begin{eqnarray*}
 \left(-\tilde\Gamma_{ik}^{j\ell} G^{km} y_m y_j \right.
  &+& \left.\Gamma_i^{j\ell}(\Gamma_j^{km}y_k y_m - g_{jm}q^m)\right)y_\ell
          -\gamma_{i\ell j}  G^{jm} y_m q^\ell + g_{ij}(G^{j\ell} p_\ell
          -\Gamma_\ell^{j m} y_m q^\ell) \\
  &=& (\Omega'' y^2)_{i}^\ell y_\ell - (\Omega' y)_{i\ell} q^\ell + g_{ij}G^{j\ell}p_\ell = L(y^3,yq,p)
\end{eqnarray*}
with $\left(\Omega'' (y_0)^2 \right)_i^{\ell}
      = -\tilde\Gamma_{ik}^{j\ell} G^{km}y_j^0 y_m^0
        + \Gamma_i^{jm} \Gamma_j^{\ell k} y^0_{m} y^0_k$.
The corresponding term for $\ddot{\tilde y}(0)$  does not vanish but the terms
cancel.
Here $L$ is a linear function of degree $3$ in its arguments with coefficients
depending on $x(t)\in\RR^d$. In the sequel we make use of the fact that we may
change the names of the indices. This way we get e.g. the equality
$\Gamma_i^{j\ell} g_{jk}q^k y_\ell = \Gamma_i^{jm} g_{j\ell}q^\ell y_m$.
Furthermore $2y_i^{(3)}(t) $ is of the form
\begin{eqnarray*}
 -\frac{\ptl \tilde\Gamma_{im}^{j\ell}}{\ptl x_n} {\dot x}^n G^{mk} y_j y_\ell y_k
  -2 \tilde\Gamma_{im}^{j\ell} G^{mk} {\dot y}_j y_\ell y_k
  - \tilde\Gamma_{im}^{j\ell}\Gamma_n^{mk} {\dot x}^n  y_j y_\ell y_k
  - \tilde\Gamma_{im}^{j\ell} G^{mk} y_j y_{\ell} {\dot y}_k
  \\
  +\tilde\Gamma_{in}^{j\ell} {\dot x}^n(\Gamma_j^{km} y_\ell y_k
  + g_{j\ell} q^\ell)y_m
  + \Gamma_i^{j\ell} (\tilde\Gamma_{jn}^{km} {\dot x}^n y_m y_k
  + 2 \Gamma_j^{km}{\dot y}_m y_k
  -\gamma_{jmk} {\dot x}^k q^m - g_{jm}{\dot q}^m)y_\ell
  \\
  +\Gamma_i^{j\ell}(\Gamma_j^{km} y_m y_k + g_{jm} q^m){\dot y}_\ell
  - \frac{\ptl\gamma_{ijm}}{\ptl x^k} G^{m\ell}{\dot x}^k  y_\ell q^j
  - \gamma_{ijm}(\Gamma_k^{m\ell}{\dot x}^k y_\ell q^j
                  +G^{m\ell}({\dot y}_\ell q^j+y_\ell {\dot q}^j))
  \\
  + \gamma_{ijk}{\dot x}^k(G^{j\ell} p_\ell - \Gamma_\ell^{ jm}y_m q^\ell)
  +g_{ij}\left(\Gamma_k^{j\ell} {\dot x}^k p_\ell + G^{j\ell} {\dot p}_\ell
  - \tilde\Gamma_{\ell k}^{jm} {\dot x}^k y_m q^\ell
  -\Gamma_\ell^{jm}({\dot y}_m q^\ell + y_m {\dot q}^\ell)\right)
\end{eqnarray*}
which constitutes a linear function $L(y^4,y^2q,q^2,yp)$ of degree $4$ in its
arguments with coefficients depending on $x(t)\in\RR^d$. For the difference
corresponding to $y^4$ and $\tilde y^4$ we give a direct method. Expanding the
tensor with respect to $x$ and $\tilde x$ and applying that we are working in Normal
coordinates we are left with an error of the order $O(t^5)$ for the difference of
these terms.
For the integrand the
estimates given in Proposition \ref{prop:2.1} apply. We find
\begin{displaymath}
  y(t)= \tilde{y}(t)+ t g(x) q_0 + \frac 12 t^2 [g(x)G(x)p^0 + (\Omega y)q_0]\ .
\end{displaymath}
As can be read from the preceding equation the higher
order terms corresponding to $\ {\ddot y}(t),\ y_i^{(3)}(t)$ adequately become
part of the error term $\delta$. In order to achieve an analogous representation
for the function $\Delta x = x - \tilde x$ we start from the equation
$\Delta x^i(t) = x^i_0 + t \dot{\Delta x}_0^i + \frac 12 t^2 \ddot{\Delta x}^i_0
  +\int_0^t (t-\tau)^3 {\Delta x^i}^{(4)}(\tau) \,d\tau$.
Here the values of the derivatives are determined by inserting into
\begin{eqnarray*}
 \ddot x^i(t)    &=& \Gamma_m^{ij} G^{m\ell} y_\ell y_j - G^{ij} \dot{y}_j= L(y^2,q)\\
  {x^i}^{(3)}(t) &=&
  -(\tilde\Gamma_{mn}^{ij}G^{m\ell} + \Gamma_m^{ij}\Gamma_n^{m\ell})G^{nk}y_k y_\ell y_j
                  + \Gamma_{m}^{ij}G^{m\ell}(\dot{y}_\ell y_j+ 2y_\ell \dot{y}_j)
                  - G^{ij} \ddot{y}_j\\
                 &=& L(y^3,yq,p)
     \\
  {x^i}^{(4)}(t)
  &=& -\left(\frac{\ptl\tilde\Gamma_{mn}^{ij}G^{m\ell}}{\ptl x^\nu}
             +\tilde\Gamma_{mn}^{ij} \tilde\Gamma^{m\ell}_{n\nu}
             + \tilde\Gamma_{m\nu}^{ij}\Gamma_n^{m\ell}
             + \Gamma_{m}^{ij}\tilde\Gamma_{n\nu}^{m\ell}
                     \right)G^{nk}G^{\nu h} y_h y_k y_\ell y_j\\
  &-&\!\!(\tilde\Gamma_{mn}^{ij}G^{m\ell}
     +\Gamma_{m}^{ij}\Gamma^{m\ell}_{n} )G^{nk} y_k(\dot{y}_\ell y_j+ 2y_\ell \dot{y}_j)
   + \Gamma_{m}^{ij}G^{m\ell}( \ddot{y}_\ell y_j +3\dot{y}_\ell\dot{y}_j +3y_\ell\ddot{y}_j)
  \\
  &-&\!\! G^{ij} y_j^{(3)} = L(y^4,y^2q,q^2,yp).
\end{eqnarray*}
Here $L(y^2,q)$, $L(y^3,yq,p)$ and $L(y^4,y^2q,q^2,yp)$ stand for different linear
functions of degree $2$, $3$ and $4$, respectively. The explicit expressions of
these functions follow by inserting the previous results on $\dot{x}$, $\ddot{x}$, and
$x^{(3)}$. For the difference of the homogenous terms $x^4$ $\tilde y^4$ we proceed
as above leaving a uniform error in time in the integrand . All together we find
\begin{equation} \label{eq:2.1d}
  x^i(t) = \tilde x^i  - \frac{t^2}6 G^{ij}g_{ij}p^0
   + \frac{t^3}6 \left(2 \frac{\ptl\Gamma_m^{ij}}{\ptl x^k} G^{kn} G^{m\ell}
   + \frac{\ptl \Gamma_m^{i\ell}}{\ptl x^k} G^{kn} G^{mj}\right)G^{jn} y_n y_\ell p_n^0\ .
\end{equation}
This concludes the proof of \ref{prop:2.2}.
\end{proof}
Using the expression for the solution of the Hamilton equation given in Proposition
\ref{prop:2.2} above we find the approximations for the corresponding derivatives w.r. to
the initial data $q_0$ and $p^0$ by plain derivation.
\begin{corollary}\label{cor:2.2i} Under the assumptions of \ref{prop:2.1}
the asumptotic behaviour of the partial derivatives of the components of the
solution to \refeq{eq:2.1} is given by
\begin{eqnarray*}
 \frac{\ptl{X}}{\ptl{p^0}} &=& \frac{1}{6} t^3 (1d + \delta), \quad
 \frac{\ptl{X}}{\ptl{q_0}} =-\frac{1}{2}t^2 (1d + \frac{t}{3}(\Omega'y^0) + \delta^2) \\ \nonumber
 \frac{\ptl{Y}}{\ptl{p^0}} &=& \frac{1}{2} t^2 (1d + \delta), \quad
 \frac{\ptl{Y}}{\ptl{q_0}} = t(1d + \frac{t}{2} (\Omega y^0) + \delta^2)  \\ \nonumber
 \frac{\ptl{P}}{\ptl{p^0}} &=& 1d + \delta, \quad
 \frac{\ptl{P}}{\ptl{q_0}} = \frac{\delta^2}{t}, \quad
 \frac{\ptl{Q}}{\ptl{p^0}}= t \cdot \delta, \quad
 \frac{\ptl{Q}}{\ptl{q_0}}= 1 + t\left(\frac{\ptl}{\ptl{x}}G(x_0)\right)y^0+\delta^2 \nonumber
\end{eqnarray*}
\end{corollary}
The bounds for the derivatives given in Proposition \ref{prop:2.2}
(and hence in its Corollary ) can be generalized to
higher order derivatives in $p^0$ and $q_0$ which will be used later. Let us
start with the $y$-component of the solution to the Hamilton equation. Instead of
the second order Taylor polynomial \refeq{eq:2.1c} we start off with a fifth order
polynomial since the second order derivative w.r. to $p^0$ needs to be carried out
explicitly. This corresponds to $y^{(5)}$ being a linear function $L$ with argument
of degree $6$. In the proof of Proposition \ref{prop:2.2} the derivatives up to
$3^{rd}$ order are calculated explicitly.
Furthermore, up to $4^{th}$ order the argument of the $n^{th}$ order derivative of
$x$ and the $(n-1)^{st}$ order derivative of $y$ are of the same degree.

For the higher order derivative we proceed in the following formal way,
in particular, our notation does not distinguish different linear
functions $L$:
\begin{eqnarray}
  y^{(5)}(t) &=&  \frac{d }{d t} L(y^4,y^2q,q^2,yp)\nonumber\\
             &=&  \frac{\ptl L}{\ptl x} (y^4,y^2q,q^2,yp)\dot{x}
                + \frac{\ptl L}{\ptl y} (y^3,yq,p)\dot{y}
                + \frac{\ptl L}{\ptl q} (y^2,q)\dot{q}
                + \frac{\ptl L}{\ptl p} (y)\dot{p}\nonumber\\
             &=& L(y^5,y^3q,yq^2,y^2p) +  L(y^5,y^3q,yq^2,y^2p,qp)\nonumber\\
             && \phantom{=L} L(y^3q,yq^2,y^2p,qp) +  L(y^3q,yq^2,y^2p) \nonumber\\
             &=& L(y^5,y^3q,yq^2,y^2p,qp)\label{eq:2.2a}
\end{eqnarray}
and analogously we find

\begin{eqnarray}
y^{(6)}(t) &=& L(y^6,y^4q,y^2q^2,q^3,y^3p,yqp,p^2)\label{eq:2.2b}\\
x^{(6)}(t) &=& L(y^5,y^3q,yq^2,y^2p,qp)\label{eq:2.2c}\\
x^{(7)}(t) &=& L(y^6,y^4q,y^2q^2,q^3,y^3p,yqp,p^2)\label{eq:2.2d}\ .
\end{eqnarray}

The algebraic structure behind the procedure is similar to
the one studied in \cite{HiL01}, \cite{NSW85}, \cite{BA89} and was developed
in the case of a regular diffusion in \cite{Kol00}.
\vskip 0.01 cm
Since the coefficients are linear combinations of
different orders of derivatives of the functions $g$ and $G$ we may conclude that
the coefficients do not vanish simultaneously.
Expanding $q$ and $p$ in a Taylor series w.r. to $t$ we find those functions $L$
which explicitly depend on $(q_0)^2$, $q_0 p^0$, and $(p^0)^2$. We are interested
in the lowest order in $t$ for which these arguments appear.

As was done in the course of the proof to Corollary \ref{cor:2.2i} for the first
order derivatives the coefficients of the second order derivatives in $q_0$ and
$ p^0$ of $X$ and $Y$ are bounded given by polynomials in $t$.
This result is summarized in the following corollary.

\begin{corollary}\label{cor:2.2ii}
Under the assumptions of \ref{prop:2.1}
the asumptotic behaviour of the second order partial derivatives of the components
of the solution to \refeq{eq:2.1} is given by
\begin{displaymath}
\begin{array}{ccc}
 \frac{\partial^2}{\partial q_0^2} X = O(t^4) &\frac{\partial^2}{\partial q_0p^0} X = O(t^5)
                                     &\frac{\partial^2}{\partial(p^0)^2} X = O(t^6)\\
 \frac{\partial^2}{\partial q_0^2} Y = O(t^3)&\frac{\partial^2}{\partial q_0p^0} Y = O(t^4)
                                     & \frac{\partial^2}{\partial(p^0)^2} Y = O(t^5) \ .
\end{array}
\end{displaymath}
\end{corollary}

The result proven in the corollary above can be summarized in the following way:
\begin{equation}\label{eq:2.3}
 \frac{\ptl (X,Y)^I}{\ptl \zeta^{I_1}\ldots\ptl \zeta^{I_k}} = O(t^{3k-I-I_1-\ldots -I_k})\ .
\end{equation}
where we have $I=0$ for $X$ and $I=1$ for $Y$ as well as $I_k=0$ for $q_0$ and $I_k=1$
for $p^0$ cf. \cite{Kol00}.

From Corollary \ref{cor:2.2ii} follows that the matrix
$\frac{\partial{(X,Y)}}{\partial{(p^0,q_0)}}$ is non degenerate hence invertible and
bounded for each $t\in(0,t_0],\ t_0<\infty$. In particular
\begin{equation}\label{eq:2.2}
 J(t,x,y,x_0,y^0) = det \frac{\partial{(X,Y)}}{\partial{(p,q)}}(t,x,y,p,q)_{
         p=p^0(t,x,y,x_0,y^0)\atop%% \\
         q=q_0(t,x,y,x_0,y^0) }
\end{equation}
does not vanish provided the conditions of Proposition 2.1 hold.
\vskip 0.2cm

\begin{theorem}\label{thm:2.3} Let the assumptions of Proposition \ref{prop:2.1} be
satisfied.
i) There exist positive real numbers $c$ and $t_0$ (uniformly in $x_0$) such that for all
$t\leq t_0$ and $|y^0|\leq \frac{c}{t}$ the mapping
$(q_0,p^0)\longmapsto (X,Y)(t,x_0,y^0,q_0,p^0)$ defined on the polydisk
$B_{\frac{c^2}{t^2}}\times B_{\frac{c^3}{t^3}}$ is a diffeomorphism onto its image.
\noindent

ii) There exists $r>0$ such that by reducing $c$ and $t_0$ if necessary, one can assume
that the image of the diffeomorphism described in i) contains a poly disc
$B_r(\tilde x)\times B_{\frac{r}{t}}(\tilde y)$.
\end{theorem}
\begin{proof}
i) Due to the asymptotic expansion for $(X,Y)$ given in Proposition \ref{prop:2.2}
the Jacobian of the mapping $(q_0,p^0)\mapsto (X,Y)(t,x_0,y^0,q_0,p^0)$
does not vanish for $t \in [0,t_0]$, $t_0>0$,  while the assumptions of
Proposition \ref{prop:2.2} hold. Hence the mapping itself is invertible. Moreover,
the mapping constitutes a diffeomorphism since we are implicitly given the
expansion for the inverse where the $t$-dependence of the poly disk compensates
for the different behaviour in time of the parameters. This finishes the proof of
the statement i).

ii) We have to prove that for arbitrary
$(x,y)\in B_r(\tilde x)\times B_{\frac{r}{t}}(\tilde y)$ there exists
$(q_0,p^0) \in B_{\frac{c^2}{t^2}}\times B_{\frac{c^3}{t^3}}$ such that
$(x,y) = (X,Y)(t,x_0,y^0,p^0,q_0)$. This is equivalent to proving the existence of
a fixed point as for the non-degenerate case treated in \cite{Kol00}.
\end{proof}
The projections on the $(X,Y)$-space of the solutions to \eqref{eq:2.1} are
called characteristics of the Hamiltonian $H$ given by \refeq{Hamilton}.
Suppose that for each initial value $(x_0,y^0) \in \mathbb{R}^{2n}$, and
$t \in [0,t_0]$, $t>0$, there exists a neighbourhood $\Omega=\Omega(t,x_0,y^0)$ of
the origin in the $(q,p)$-space $\mathbb{R}^{2n}$ such that the mapping
$(q_0,p^0)\mapsto (X,Y)(t,x_0,y^0,q_0,p^0)$ is a diffeomorphism from $\Omega$ onto
its image containing a neighbourhood $D(x_0,y^0)$, then the family
% 2.4
\[ \Gamma(x_0,y^0)=\{(X,Y)(t,x_0,y^0,q_0,p^0)| (q_0,p^0)\in \Omega(x_0,y^0)\}, 0\leq t\leq t_0
\]
of solutions  is called a field of characteristics starting at $(x_0,y^0)$. Due to
Corollary \ref{cor:2.2i} and the subsequent arguments this means that under the
assumptions of Proposition \ref{prop:2.1} there exists a field of characteristics
for each initial point $(x_0,y^0)$ i.e. there exists a smooth function
\[ (Q_0,P^0)(t,x,y,x_0,y^0): (0,t_0] \times D(x_0,y^0) \rightarrow \Omega
\]
such that
\begin{equation}\label{eq:2.4}
       (X,Y) \left( t,x_0,y^0,(P^0,Q_0)(t,x,y,x_0,y^0) \right) = (x,y)\ .
\end{equation}
Recalling that there exists a field of characteristics for any two points
$(x_1,y_1)$ and $(x_2,y_2)$ joint by a characteristic curve of
the Hamiltonian $H$  we
introduce the function
\begin{equation}
  \sigma(t,x_0,y^0,q_0,p^0)
     = \int_0^t\left((P,Q),(\dot X,\dot Y)\right)(\tau)-H(X,Y,P,Q)(\tau) d\tau\ ,
\end{equation}
where the components $X,Y,P,Q$ of the solution to the Hamiltonian system
corresponding to $H$ are being evaluated at an intermediate time
\mbox{$\tau, \,\, 0\leq \tau \leq t$}. Due to Proposition \ref{prop:2.1}
This function is diffeomorphic equivalent to
\begin{equation} \label{eq:2.5}
   S(t,x,y,x_0,y^0) = \int_0^t (Q,P)(\tau) d\gamma
                      - \int_0^t H(\gamma,P,Q)(\tau) d\tau\ ,
\end{equation}
which is called the \emph{two point function}.

\begin{proposition}\label{prop:2.5}
Under the assumptions of Proposition~\ref{prop:2.1} and Theorem~\ref{thm:2.3} the
function $S(t,x,y,x_0,y^0)$ satisfies the Hamilton Jacobi equation
\begin{displaymath}
  \frac{\partial S}{\partial t}
  + H(x,y,\frac{\partial S}{\partial x},\frac{\partial S}{\partial y}) = 0
\end{displaymath}
as a function of $t \in (0,t_0]$ and $x,y \in D(x_0,y^0)$ with $D(x_0,y^0)$ as
above. Moreover, we have
\begin{eqnarray*}
  \frac{\ptl S}{\ptl x}(x,y)   &=& p(t,x,y), \quad \frac{\ptl S}{\ptl y}(t,x,y)=q(t,x,y).\nonumber\\
 \frac{\ptl S}{\ptl x_0}(x,y), &=& -p^0(t,x,y),
         \quad \frac{\ptl S}{\ptl y^0}(t,x,y)= - q_0(t,x,y)\nonumber
\end{eqnarray*}
\end{proposition}
Idea of proof: The partial derivatives are considered first. Since the inverses
$\frac{\ptl X}{\ptl q_0}^{-1} $ and $ \frac{\ptl Y}{\ptl p^0}^{-1}$ exist due
to Corollary \ref{cor:2.2i} the functions $S$ and $\sigma$ are identical up to the
diffeomorphism given in Theorem \ref{thm:2.3}. Now the derivatives follow by direct
calculation. That means that from this point the proof follows the
line of arguments for the nondegenerate variational principle, see e.g. \cite{Kol00}.

\begin{proposition}\label{prop:2.3}
Let the assumptions of Proposition~\ref{prop:2.1} be satisfied, then the following
representations hold:
\begin{displaymath}
\begin{array}{rclrcl}
\frac{\ptl^2{S}}{\partial{x^2}} &=&
    \frac{\partial{P}}{\partial{p^0}}\left(\frac{\partial{X}}{\partial{p^0}}\right)^{-1}, &\quad
\frac{\partial^2{S}}{\partial{y^2}} &=&
    \frac{\partial{Q}}{\partial{q_0}}\left(\frac{\partial{Y}}{\partial{q_0}}\right)^{-1}, \\
\frac{\partial^2{S}}{\partial{x}\partial{y}} &=&
\frac{\partial{Q}}{\partial{p^0}}\left(\frac{\partial{X}}{\partial{p^0}}\right)^{-1}, &\quad
\frac{\partial^2{S}}{\partial{y}\partial{x}} &=&
\frac{\partial{P}}{\partial{q_0}}\left(\frac{\partial{Y}}{\partial{q_0}}\right)^{-1},\\
\end{array}
\end{displaymath}

where the functions $S,X,P,Q$ depend on the same tuple $(t,x,y,q_0,p^0)$ of variables.
\end{proposition}
This is achieved by iteration from Proposition \ref{prop:2.5} and Corollary \ref{cor:2.2i}
which implies Theorem \ref{thm:2.3}.
\begin{remark}\label{rem:3.1} Due to the symmetry of the results in Corollary \ref{cor:2.2i} and
Proposition \ref{prop:2.5}
an analogue to Proposition \ref{prop:2.3} in the variables $x_0,\, y^0$ holds.
For the lowest order approximation the following explicit expressions exist for
each component:
\begin{eqnarray}\label{invcovariance}
\frac{\partial^2{S}}{\partial^2{\{x,y\}}}=
\left(\!\!
  \begin{array}{cc}
    \frac{12}{t^3}(1+\delta) & \frac{-6}{t^2}(1+\delta) \\
    \frac{-6}{t^2}(1+\delta) & \frac{4}{t}(1+\delta) \\
  \end{array}
\!\!\right)
\quad
\mbox{and}
\quad
\frac{\partial^2{S}}{\partial^2{\{x_0,y^0\}}}=
\left(\!\!
  \begin{array}{cc}
    \frac{12}{t^3}(1+\delta) & \frac{6}{t^2}(1+\delta) \\
    \frac{6}{t^2}(1+\delta) & \frac{4}{t}(1+\delta) \\
  \end{array}
\!\!\right) .
\end{eqnarray}
\end{remark}
For the Hamiltonian $H(\xi, \eta)$ let us define the functional
\begin{equation}\label{eq:2.6}
  I_t(\xi) = \int_0^t \max_{\eta}
             \left\{ (\eta, \dot \xi)-H  (\xi, \eta)\right\}(\tau) d\tau
\end{equation}
The maximum is taken over all piecewise smooth curves $\xi$ with given endpoints
$\xi(0)=\xi_0$ and $\xi(t)=\varkappa$.
\begin{proposition}\label{prop:2.2ii}
For the Hamiltonian given in \eqref{Hamilton} the characteristic curve \\
$(X,Y)(t,x_0,y^0,q_0,p^0)$
provides a minimum for the functional (\ref{eq:2.6}) over all curves in
$D(x_0,y^0)=B_{\frac{c^2}{t^2}}(q_0) \times B_{\frac{c^3}{t^3}}(p^0)$.
Furthermore, $S(t,x,y,x_0,y^0)$, of \eqref{eq:2.5}, constitutes the unique
minimal value of the functional \refeq{eq:2.6}.
\end{proposition}
\begin{proof}
For the Hamiltonian given in (\ref{Hamilton}) the Weierstra{\ss} condition is
satisfied, i.e. for all $t>0$ and all $\xi,\ \eta\in\RR^d$:
$W(x(t),y(t),\xi,\eta,p(t),q(t)) \ge 0$
 for arbitrary solutions $(x(t),y(t),q(t),p(t))$ of \refeq{eq:2.1} where
the Weierstra{\ss} function is given by
\begin{displaymath}
  W(x,y,\xi,\eta,p,q) = \Delta H  - L(H)
\end{displaymath}
with $H$ as defined in \refeq{Hamilton}, $\Delta H:= H(x,y,q_1,p_1)-H(x,y,q_2,p_2)$ and $L(H)=$
$\left(\Delta q,\frac{\ptl H}{\ptl q}\right)+\left(\Delta p,\frac{\ptl H}{\ptl p}\right)$
$:= \left(q_1-q_2,\frac{\ptl H}{\ptl q}\right)+\left(p_1-p_2,\frac{\ptl H}{\ptl p}\right)$.
Then we find
\begin{eqnarray*}
 \Delta H &=& \frac 12 \left(g(x)q_1,q_1 \right)- \left(g(x)q_2,q_2 \right)-
 \left(G(x)y,\Delta p\right)-\left(\frac{\ptl}{\ptl x}(G(x)y,y),\Delta q \right)\\
      L(H)&=& \left(q_1-q_2,g(x)q_2 \right)
      + \frac 12\left(\frac{\ptl}{\ptl x}(G(x)y,y),\Delta q \right) -(G(x)y,\Delta p)\ .
\end{eqnarray*}
Combining the two equations reveals
\begin{displaymath}
 W= \frac 12 (g(x)q_1,q_1 ) - (q_1,g(x)q_2)+ \frac 12 (g(x)q_2,q_2)
   = \frac 12 (q_1-q_2,g(x)(q_1-q_2)) + O(1) \ge 0\, ,
\end{displaymath}
the positivity holding since the matrix $g$ is positive definite. As a consequence the Hamilton function treated
in this paper satisfies the Weierstrass condition. Inserting definitions and Theorem
\ref{thm:2.3} we have $I_t\ge S$ c.f. \refeq{eq:2.5}. Hence the minimum for the functional
will be unique for arbitrarily chosen characteristics in
$B_{\frac{c^2}{t^2}}(q_0) \times B_{\frac{c^3}{t^3}}(p^0)$, which finishes the proof
of Proposition~\ref{prop:2.2ii}.
\end{proof}
Combining the result of Theorem \ref{thm:2.3} with the
uniqueness result obtained in Proposition~\ref{prop:2.2ii} gives us the equivalence of
solutions to the Hamilton equation and Hamilton Jacobi equation.
\vskip 0.1 cm
We conclude this section by showing that the two-point function $S$ is convex on $D(x_0,y^0)$.
For arbitrary square matrices $a,b,c,d$ the general following formula holds
\[
  det \left(
      \begin{array}{cc}
        a & b \\
        c & d
      \end{array} \right)
      = det~d \quad det~(a-bd^{-1}c)\ .
\]

Using the asymptotics given in Remark \ref{rem:3.1} as a starting point, applying
the above result and integrating twice w.r. to $x$ and $y$ we see that the two
point function can be uniformly approximated by a quadratic function in $x$ and $y$.
In particular the two point function $S$ can be approximated as follows
\begin{eqnarray}\label{eq:2.7}
  S(x,y,x_0,y_0) = \frac{6}{t^3}(1+\delta)(x-x_0)^2
             + \frac{6}{t^2} (1+\delta) \left( (x-x_0), y+y^0 \right)
      \nonumber \\
      +  \frac{2}{t} (1+\delta)[(y,y) + (y,y^0) + (y^0,y^0)].
\end{eqnarray}
This leads to a Gaussian approximation of the solution to \refeq{Greensfunktion}.

\section{Uniform Estimates for the Heat Kernel}
In this section we shall construct an asymptotic expansion for the fundamental
solution (Greens function) of the degenerate diffusion \eqref{Greensfunktion}.
To this end the well known WKB-type method is modified in such a way
that the variables representing derivatives i time of order zero to three
are being transformed in order to lift them to the same time scale. Hereby closed
algebraic equations for the coefficients of the expansion are achieved, see also
the construction in \cite{Kol00} for the origin of our construction.

We shall discuss two types of asymptotics for the solution to
\eqref{Greensfunktion} and \eqref{Hamilton}, namely, small time asymptotics when
$t \to 0$ and $h$ is fixed ($h=1$) and small diffusion asymptotics when $t$ is
fixed and $h \to 0$.

We look for the asymptotic fundamental solution to \refeq{eq:1} using a
WKB Ansatz
\begin{equation}\label{eq:3.1}
  u^0(t,x,y,x_0,y^0,h) = C(h)\phi(t,x,y,x_0,y^0) \exp[-S(t,x,y,x_0,y^0)/h]\ ,
\end{equation}
with real functions $\phi$ and $S$ called the amplitude and the phase, respectively.
Moreover, $C$ is a function only depending on $h$ alone while $\phi$ and $S$ are
independent of $h$ having the arguments $x,y,x_0,y^0\in \RR^d.$
This particular form of the asymptotic fundamental solution applies due to the
positivity of the resolving diffusion operator. It is specific to linear second
order pseudo differential equations describing Markov processes.

This asymptotic ansatz is also called exponential.
Hence one should require that the solution $u$ of
\eqref{Greensfunktion} differs from $u^0$ by a multiplicative error:
\begin{equation}\label{eq:3.3}
  u(t,x,y,x_0,y^0,h) = C(h) \phi(t,x,y,x_0,y^0) \exp[-S(t,x,y,x_0,y^0)/h](1 + O(h)).
\end{equation}
For \eqref{eq:3.3} only the term with minimal entropy of the sum survives at each
point.
Therefore, for \eqref{eq:3} the common superposition principle transforms to the
idempotent superposition principle $(S_1,S_2) \mapsto min(S_1,S_2)$ at the level of
actions. For a detailed discussion of this idempotent superposition principle
and its applications see \cite{KoM94}, \cite{KoM97}.
Inserting $u^0$ \eqref{eq:3.1} into \eqref{Greensfunktion} reveals

\begin{eqnarray}\label{eq:3.4}
  h\left(\frac{\partial \phi}{\partial t} -\frac 1h \phi \frac{\partial S}{\partial t}\right)
  = \frac{h^2}2 tr~ g(x)\left(\frac{\partial^2 \phi}{\partial y^2}
   -\frac\phi h \frac{\partial^2 S}{\partial y^2}\right)
   + h\left(G(x)y,\frac{\partial \phi}{\partial x}
   -\frac 1h \phi \frac{\partial S}{\partial x}\right) \\
   + \frac{-h}2 \left(\frac{\partial}{\partial x}(G(x)y,y),\frac{\partial \phi}{\partial y}
   - \frac 1h \phi \frac{\partial S}{\partial y}\right)
   + \frac 12 \left(g(x)\frac{\partial S}{\partial y}\right)\phi
   - h\left(g(x)\frac{\partial S}{\partial y}, \frac{\partial \phi}{\partial y}\right) \nonumber
\end{eqnarray}

The terms in the above equation may be classified according to the degree of $h$
which ranges from $0$ to $2$. Comparing coefficients we obtain  the following two
equations: The first Hamilton-Jacobi equation ($h$ to the order $0$):
\begin{equation}\label{Hamilton-Jacobi}
  \frac{\partial S}{\partial t} + \frac 12
\left(g(x)\frac{\partial S}{\partial y},\frac{\partial S}{\partial y}\right)
- \left(G(x)y, \frac{\partial S}{\partial y}\right)
+ \frac 12 \left(\frac{\partial }{\partial x} (G(x)y,y), \frac{\partial S}{\partial y}\right) = 0
\end{equation}
which in terms of the Hamiltonian function \eqref{Hamilton} gives
\begin{displaymath}
  \frac{\partial}{\partial t}S +
  H\left(x,y, \frac{\partial S}{\partial x},\frac{\partial S}{\partial y}\right) = 0 \quad .
\end{displaymath}
This way the system given by the PDE \refeq{eq:1} is linked with the Hamilton-Jacobi
equation studied in Section $2$.
The second equation is  the so called transport equation (the coefficient of $h^1$):
\begin{equation}\label{Transport}
  \frac{\partial \phi}{\partial t}
  + \left(\!g(x)\frac{\partial S}{\partial y},\frac{\partial \phi}{\partial y}\right)
  + \frac 12 tr\!\left(\!g(x)\frac{\partial^2 S}{\partial y^2}\right)
  -~\left(G(x)y, \frac{\partial \phi}{\partial x}\right)
  +\frac 12\left(\!\left(\frac{\partial G(x)}{\partial x}y,y\right), \frac{\partial \phi}{\partial y}\right)
  = 0.
\end{equation}
Going the opposite way by choosing $S$ and $\phi$ which satisfy \eqref{Hamilton-Jacobi}
and \eqref{Transport} and inserting then into the r.h.s of \eqref{Greensfunktion} we
define a function $u^0$ which satisfies equation \eqref{Greensfunktion} up to an error of order
$h^2$. More
precisely we obtain the following inhomogeneous equation
\begin{equation}\label{uNULL}
 h\frac{\partial u^0}{\partial t}
 \! -\! H\!\left(\!x,y,-h\!\frac{\partial}{\partial x},-h\!\frac{\partial}{\partial y}\!\right)\! u^0
 \!\! =\!\!
  -\frac{h^2}2 C(h)F_0(t,x,y,x_0,y^0)
\end{equation}
$F_0(t,x,y,x_0,y^0)\equiv
                   tr\!\left(\!g(x)\frac{\partial^2 \phi}{\partial y^2}\!\right)
    \exp(-\frac{S(t,x,y)}{h})$
i.e. the right hand side defines a function $F_0(t,x,y,x_0,y^0)$. If we would concentrate
on a Cauchy problem instead with a smooth initial function then we would have to solve an
appropriate Cauchy problem for the Hamilton-Jacobi equation.\\
We now proceed to derive the precise asymptotic behaviour of the solution to the
parabolic equation (\ref{Greensfunktion}) for small times (respectively in a small
neighbourhood of the starting point) we shall provide  asymptotic expansions in
the time parameter $t$ for the two point function $S$ and the amplitude $\phi$
given in \refeq{eq:3.1}. Studying the inhomogeneous system \refeq{uNULL} has the
advantage that the functions $S$ and $\phi$ satisfy the Hamilton-Jacobi- \refeq{Hamilton-Jacobi}
respectively the transport equation \refeq{Transport}.

In Section 3 of this paper we have seen that the two point function $S$ in
\ref{eq:2.5} is smooth and satisfies the Hamilton-Jacobi equation almost everywhere.
Moreover, it was shown that the solution of the Hamilton-Jacobi equation and the
Hamilton equations are equivalent thus generalizing the results given in
\cite{Kol00} for non degenerate systems to a special case of a regular Hamiltonian
of rank one. Due to the results in Section 3 we only use the expansions of the
solution to the Hamiltonian system in the sequel.

In a fist step we neglect the inhomogeneity in \refeq{uNULL}.
Starting with the phase respectively the two point function we are looking for a
representation of the form
\begin{equation}\label{phase}
  S(t,x,y) = \sum_{i=-k}^\infty S_{i}(x,y)t^{i}   \quad .
\end{equation}
Traditionally this is achieved by introducing the asymptotic expansion for the
corresponding Hamiltonian system \refeq{eq:2.1} into the Hamilton-Jacobi equation
at the point $(t,x,y)$ and subsequently solving a closed system  of algebraic equations in
$x$ and $y\in\mathbb{R}^d$. This fails in the case of our degenerate system. In
this paper we proceed by introducing a change of variables which
first shifts the origin into the solution $(\tilde{x},\tilde{y},\tilde{p},\tilde{q})$
of the Hamilton equation \eqref{eq:2.1} starting at $(x_0,y_0,0,0)$
and then adjusts the time scales of the variables
on the manifold and the tangent space. We define a function $\sum$ on
$\RR_+\times\RR^d\times\RR^d$ by
\begin{equation}\label{eq:3.9}
  \sum(t,\xi,y) = S(t,t(x+\tilde{x}(t)),y + \tilde y ,x_0, y^0)\ .
\end{equation}
We remark that the function
$\sigma(t,x,y) = S(t,x+\tilde{x}(t),y + \tilde y ,x_0, y^0)$ coincides with the one
given in \refeq{eq:2.5}. There holds
$  \frac{\partial \Sigma}{\partial \xi} =
  t\frac{\partial \sigma}{\partial x}\quad\mbox{and}\quad \frac{\partial \Sigma}{\partial t} =
  \frac{\partial \sigma}{\partial t} + x \frac{\partial \sigma}{\partial x}\ .
$
In this way a closed system of algebraic equations for the
coefficient functions $S_{i}$ of an expansion of the phase function in powers of
$x$ and $y$ will be derived.

The equation \refeq{Hamilton-Jacobi} becomes
\begin{eqnarray}\label{eq:3.11}\lefteqn{
  0 =\frac{\partial \Sigma}{\partial t} - \frac{\xi + G(t\xi + \tilde{x})(y+\tilde{y})
 - G(\tilde{x})\tilde{y}}t \frac{\partial \Sigma}{\partial \xi}
 - \left(g(\tilde{x})\tilde{q}, \frac{\partial \Sigma}{\partial y}\right)}
  \\ \nonumber
 &+&\!\!\frac 12\left[\!\left(\!\frac{\partial G(t\xi+ \tilde{x})}
                   {\partial x}(y+\tilde{y}),(y+\tilde{y})\!\right)
 - \left(\frac{\partial G(\tilde{x})}{\partial x}\tilde{y},\tilde{y}\right)\!\right]
  \frac{\partial \Sigma}{\partial y}
 + \frac 12 \left(\!g(t\xi + x)
    \frac{\partial \Sigma}{\partial y},\frac{\partial \Sigma}{\partial y}\!\right)
\end{eqnarray}
\begin{proposition}\label{prop:3.1}
For the transformed two-point function $\Sigma$ as defined in \eqref{eq:3.9} there
exists a unique asymptotic expansion in powers of $t > 0$ of the form
\begin{displaymath}
  \Sigma = \frac{\Sigma_{-1}}t + \Sigma_0 + t\Sigma_1 + t^2\Sigma_2 + ...
\end{displaymath}
such that $\Sigma_1$ and $\Sigma_0$ vanish at the origin, $\Sigma_1$ is strictly
convex in a neighbourhood of the origin and all $\Sigma_j,\ j\ge 2$, are regular power
series in $(\xi,y)$. Moreover, the remainder of arbitrary order can be estimated
by means of induction using Proposition \ref{prop:2.2} as initialization.
\end{proposition}
{\bf Proof}
As in Section 3 let us replace the coefficient functions $g,G$ in equation
\ref{eq:3.11}
by their $n^{th}$ order Taylor series and let us also insert the $n'^{th}$ order
approximation of $(\tilde{x},\tilde{y},\tilde{p},\tilde{q})$. Then the asymptotic
result given in Corollary \ref{cor:2.2ii} gives that the combined error due to this
replacement is of the order $0(t^{n+1})$.

This way \eqref{eq:3.11} becomes in coordinate form (normal coordinates) for $n=1$
up to an error of order $0(t^2)$:
\begin{eqnarray}\label{eq:3.13}
\lefteqn{\frac{\partial \Sigma}{\partial t}
  - \frac{(\xi + y)_i + \frac {t^2}2 g_{ij}^{k\ell}[(\xi_k - y_k^0)(\xi_\ell - y_\ell^0)(y_j + y_j^0) - y_k^0y_\ell^0y_j^0] + O(t^3)}t \frac{\partial \Sigma}{\partial \xi_i}}\nonumber\\
& &-~\frac 12[t g_{ij}^{k\ell}[(\xi_\ell - y_\ell^0)(y_i + y_i^0)(y_j + y_j^0) - y_i^0y_j^0y_\ell^0] + O(t^2)]\frac{\partial \Sigma}{\partial y_k}  \\
& &+~\frac 12 (1 + \frac{t^2}2 g_{ij}^{k\ell}(\xi_{k} - y_{k}^0)(\xi_{\ell k} - y_{\ell}^0) + O(t^3))\frac{\partial
\Sigma}{\partial y_{i}}\frac{\partial \Sigma}{\partial y_{j}} = 0\ .\nonumber
\end{eqnarray}
We now rearrange the l.h. side of this equation according to the order in $t$.
For the lowest orders in $t$ we arrive at the subsequent set of equations. In particular
for the lowest order in t, i.e. $t^{-2}$, equation (\ref{eq:3.13}) reduces to the following
first order homogeneous partial differential equation with linear coefficient functions:
\begin{equation}\label{eq:3.14}
  -{\Sigma}_{-1}-(y+\xi)\frac{\partial{\Sigma}_{-1}}{\partial \xi}
  + \frac 12(g_0 \frac{\partial{\Sigma}_{-1}}{\partial y}, \frac{\partial {\Sigma}_{-1}}{\partial y}) = 0\ .
\end{equation}
This linear first order partial differential equation is of the following general form
\begin{equation}\label{eq:3.11a}
  \lambda u + (Ax,\nabla u) = p(x)
\end{equation}
where $\lambda \in \mathbb R^+,~ x\in\mathbb R^m,~ p(x)$ is a homogeneous polynomial, and $A$
is an $m\times m$ matrix with strictly positive eigenvalues
$a_1\leq a_2\leq \ldots \leq a_m$. The solution constitutes a polynomial of a particular
form given in the following lemma. The proof goes by inserting and direct calculations.
\begin{lemma}\label{lemma:3.1}
i) Let p be a homogeneous polynomial of degree $q$. Then the solution of \eqref{eq:3.11a}
exists, moreover it constitutes a polynomial of degree $q$ with coefficients of the form
\begin{displaymath}
  \frac {\ptl^q u}{\ptl x_{i_1}\ldots \ptl x_{i_q}}
  = \frac 1{a_{j_1}+\ldots +a_{j_q}} (C^{-1}_{j_1i_1})\cdots (C^{-1}_{j_qi_q})
    \frac {\ptl^q p}{\ptl x_{\ell_1}\ldots \ptl x_{\ell_q}} C_{j_1i_1}\cdots C_{j_qi_q}\ .
\end{displaymath}
Let $p=\sum_{q=0}^m p_q$ be a sum of homogeneous polynomials of degree $q$ (in
case $m=\infty$ the sum has to be absolutely convergent in a ball). Then the
analytic solution of \eqref{eq:3.11a} exists and is given by the sum
$\sum_{q=1}^m u_q$ of solutions $u_q$ corresponding to the inhomogeneity $p_q$. In
case $m=\infty$ the domains of convergence of $p$ and the solution $u$ to
\refeq{eq:3.11a} coincide.
\end{lemma}
For a proof see \cite{Kol00}. We now proceed to calculate the coefficient functions
$\Sigma_i$, $-1\le i$.
\begin{lemma}\label{lemma:3.2}
Under the additional assumption ${\Sigma}_{-1}(0,0) = 0$ and that ${\Sigma}_{-1}$ is strictly
convex at the origin the formal power series for ${\Sigma}_{-1}$ in powers of $\xi$
and $y$ reduces to the quadratic form $6(\xi,\xi) + 6(\xi,y) + 2(y,y)$.
\end{lemma}
\begin{proof}
Applying the assumption that ${\Sigma}_{-1}$ has an asymptotic expansion
\begin{displaymath}
  {\Sigma}_{-1} = \sum_{k\geq 0} \zeta_k
\end{displaymath}
where $\zeta_k$ is a polynomial of order k in $\xi$ and $y$ and using the assumption
${\Sigma}_{-1}(0,0)=0$ it immediately follows that $\zeta_0 = 0$. In local coordinates
$g_0=1$.
The first order term is of the form $\zeta_1= K \xi + Ly$ where $K,L\in \RR$.
By inserting into \eqref{eq:3.14} and using that $\zeta_0=0$ we get the equation
\begin{displaymath}
 -K\xi - L y -(y+\xi) K + \frac 12 L^2 =0\ .
\end{displaymath}
At the origin we find $L^2 = 0$. Differentiating with respect to $y$
reveals $-K-L=0$ and hence $\zeta_1$ vanishes altogether, i.e  $\zeta_1\equiv 0$.
Assuming for the quadratic part $\zeta_2$ the representation
\begin{equation}\label{eq:3.15}
  -\frac 12 (A\xi, \xi) - (B\xi, y) + \frac 12 (Cy, y)
\end{equation}
we get the equations
\begin{displaymath}
  A=\frac 12 B^tB\quad 2B+A = CB\quad \frac 12 C+B = \frac 12 C^2.
\end{displaymath}
Solving this algebraic system of equations we find, since \eqref{eq:3.15} should be
positive definite, $C = 4\ id,\quad B=6\ id, ~\mbox{and}\quad A=12\ id$, i.e.
\eqref{eq:3.15} is of the form
\begin{equation}\label{eq:3.16}
  6(\xi, \xi) + 6(\xi, y) + 2(y,y)\quad .
\end{equation}
Finally all $\zeta_k,\, k\geq 3$, vanish, in particular for $k=3$ the terms of
order 3 in $\xi,y$ make up the equation
\begin{equation}\label{eq:3.17}
  0=-\zeta_3-\left((y+\xi),\frac{\partial \zeta_3}{\partial \xi}\right)
    + \left(-\frac{\partial \zeta_2}{\partial \xi},-\frac{\partial \zeta_3}{\partial y}\right) =
    \zeta_3 + \left((y+\xi) ,\frac{\partial \zeta_3}{\partial \xi}
         - (6\xi+4y),\frac{\partial \zeta_3}{\partial \xi}\right)\ .
\end{equation}
where \eqref{eq:3.16} has been inserted and the right hand side has been multiplied
by $-1$.
Due to Lemma \ref{lemma:3.1} this implies that $\zeta_3 = 0$. Inserting  $\zeta_k$,
$k\geq 3$, into (\ref{eq:3.14}) reproduces
\eqref{eq:3.17}, hence $\zeta_k$ vanish f for $k\geq 3$, which finishes the proof
of the lemma.
\end{proof}

For the term $\Sigma_0$ corresponding to the order $t^0$ in the expansion of
$\Sigma$ we have the equation
\begin{displaymath}
  -(y+\xi)\frac{\partial {\Sigma}_0}{\partial \xi}
          + (6\xi + 4y) \frac{\partial {\Sigma}_0}{\partial y} = 0\ ,
\end{displaymath}
the solution of which vanishes due to Lemma \ref{lemma:3.1} since the equation is
homogeneous. For the coefficients of $t$ in \refeq{eq:3.13} we obtain the equation
\begin{equation}\label{eq:3.17a}
  {\Sigma}_1 - (y+\xi)\frac{\partial \Sigma_1}{\partial \xi}
         + (6\xi + 4y)\frac{\partial \Sigma_1}{\partial y} = F_1(\xi, y)\ ,
\end{equation}
where $F_1(\xi,y)$ is the sum of the homogeneous polynomials $F_{11}, F_{12}$ and $F_{13}$ of
order 2,3 and 4, respectively, in $\xi$ and $y$. In particular we have
\begin{eqnarray*}
F_{11} &=& g_{ij}^{kl}[(12\xi_i\xi_k-4y_iy_k)y_j^0y_l^0-(18\xi_iy_j+7y_iy_j+9\xi_i\xi_j)y_k^0y_l^0+(3\xi_k\xi_l+2\xi_ky_l)y_i^0y_j^0]\\
F_{12} &=& g_{ij}^{kl}[(-6\xi_i\xi_k\xi_l+3\xi_k\xi_lly_i+4\xi_ky_iy_l)y_j^0\\
     && \qquad +(36\xi_i\xi_ky_j+11\xi_ky_iy_j-2y_iy_jy_k+18\xi_i\xi_j\xi_k)y_l^0]\\
F_{13} &=& g_{ij}^{kl}[2\xi_ky_iy_jy_l-4\xi_k\xi_ly_iy_j-18\xi_i\xi_k\xi_ly_j-9\xi_i\xi_j\xi_k\xi_l]\ .
\end{eqnarray*}
Due to Lemma \ref{lemma:3.1} the form of the solution $\Sigma_1$ is determined to
be a polynomial of degree 4 in $\xi, y$ and $y^0$. The expansion of $\Sigma$ in
powers of $t$ (having expanded in $x$ and
$y$ before) is explicitly given by

\begin{displaymath}
\begin{array}{l}
\Sigma=  {\frac {6\,{\xi_{{i}}}^{2}+6\,\xi_{{i}}y_{{i}}+2\,{y_{{i}}}^{2}}{t}}+
t \left(-{\frac {14}{15}}\,g_{i,j}^{k,l}y_{0,k}y_{0,l}y_{{i}}y_{{j}}
-{\frac {59}{20}}\,g_{i,j}^{,k,l}y_{0,k}y_{0,l}\xi_{{i}}y_{{j}}
+\frac{3}{5}\, g_{i,j}^{k,l}y_{0,k}y_{0,l}\xi_{{i}}\xi_{{j}}\right.\\[2mm]
\left.\mbox{}-\frac{2}{15}\,g_{i,j}^{k,l}y_{0,j}y_{0,l}y_{{i}}y_{{k}}
+\frac{7}{5}\,g_{i,j}^{k,l}y_{0,j}y_{0,l}\xi_{{i}}y_{{k}}
+ \frac{24}{5}\,g_{i,j}^{k,l}y_{0,j}y_{0,l}\xi_{{i}}\xi_{{k}}
+\frac{1}{15}\,g_{i,j}^{k,l}y_{0,i}y_{0,j}y_{{l}}y_{{k}}\right.\\[2mm]
\left.\mbox{}+{\frac {11}{20}}\,g_{i,j}^{k,l}y_{0,i}y_{0,j}\xi_{{k}}y_{{l}}
+\frac{3}{5}\,g_{i,j}^{k,l}y_{0,i}y_{0,j}\xi_{{k}}\xi_{{l}}+
{\frac {7}{10}}\,g_{i,j}^{k,l}y_{0,l}\xi_{{i}}y_{{k}}
y_{{j}}-{\frac {9}{10}}\,g_{i,j}^{k,l}y_{0,l}\xi_{{i}}\xi_{{j}}y_{{k}}\right.\\[2mm]
\left.\mbox{}-\frac{3}{2}\,g_{i,j}^{k,l}y_{0,l}\xi_{{i}}\xi_{{j}}\xi_{{k}}
-\frac{1}{5}\,g_{i,j}^{k,l}y_{0,j}\xi_{{i}}y_{{l}}y_{{k}}
-\frac{6}{5}\,g_{i,j}^{k,l}y_{0,j}\xi_{{i}}\xi_{{k}}y_{{l}}-3\,g_{i,j}^{k,l}y_{0,j}\xi_{{i}}\xi_{{k}}\xi_{{l}}
-{\frac {3}{56}}\,g_{i,j}^{k,l}\xi_{{i}}y_{{l}}y_{{k}}y_{{j}}\right.\\[2mm]
\left.\mbox{}+\frac{1}{7}\,g_{i,j}^{k,l}\xi_{{i}}\xi_{{j}}y_{{l}}y_{{k}}
+\frac{5}{7}\,g_{i,j}^{k,l}\xi_{{i}}\xi_{{j}}
\xi_{{k}}y_{{l}}+g_{i,j}^{k,l}\xi_{{i}}\xi_{{j}}\xi_{{k}}\xi_{{l}}\right)
+O(t^2) \ .
\end{array}
\end{displaymath}
The calculations for higher orders go by analogy although they become substantially wilder.
The terms $\Sigma_j,~j>1$, turn out to be homogeneous polynomials in $\xi, y, y^0$ of
degree $j+3$
which can be estimated analogously to the term $\Sigma_1$.

In normal coordinates in particular the
term $\Sigma_1$ is given by $2.5 g^{kl}_{ij} y^0_i y^0_j y^0_k y^0_l$ due to the symmetry
of $g^{kl}_{ij}$ which is why we omit it when transforming back.
The inverse transformation $S(t,x,y,y^0) = \Sigma(t,\frac{x-\tilde{x}(t)}t , y-\tilde{y}(t))$ has
been found explicitly for an approximation of order 1 using MAPLE, a programme for symbolic and
numerical algebra, as being

\begin{eqnarray*}
\lefteqn{S_i(t,x,y,x_0,y^0)=t^{-1}\left(-6\,{\frac {{\it x}_{{0,i}}y^0_{{i}}}{t}}
-6\,{\frac {{\it x}_{{0,i}}{\it y}_{{i}}}{t}}
+6\,{\frac {x_{{i}}{\it y}_{{i}}}{t}}
+6\,{\frac {x_{{i}}y^0_{{i}}}{t}}\right.}\\
&-&\left.12\,{\frac {x_{{i}}{\it x}_{{0,i}}}{{t}^{2}}}
+6\,{\frac {{{\it x}_{{0,i}}}^{2}}{{t}^{2}}}
+6\,{\frac {{x_{{i}}}^{2}}{{t}^{2}}}
+2\,{y^0_{{i}}}^{2}
+2\,{{\it y}_{{i}}}^{2}
+2\,y^0_{{i}}{\it y}_{{i}}+O(t)\right)\ .
\end{eqnarray*}
This is consistent with the estimate found in (\ref{invcovariance}).
On the diagonal this expression is of the form $S_i= 6 y_i^2 t^{-1}$ which
obviously is not vanishing and hence reveals an exponential decay for the kernel
on the diagonal.

For $x_0=0$ this reads:
\begin{eqnarray*}
S(t,x,y,y^0) &=& \frac 6{t^3}[1+O(|x|^2) + O(|y-y^0|^2)]|x|^2\\
  & & +~\frac 6{t^2}\langle x,y+y^0+O(|x|^2)\rangle
         %\quad \rightsquigarrow x =(|x|^2)???
         \\
  & & +~\frac 2{t}[|y|^2 + \langle y^0,y \rangle + |y^0|^2 + O(|x|^2) + O(|y-y^0|^2)]\ .
\end{eqnarray*}

Having completed the first step in giving an approximate solution of exponential
WKB type \refeq{eq:3.1}, namely, having derived an asymptotic expansion for the solution $S$ to
\refeq{Hamilton-Jacobi} we now turn to the transport equation \refeq{Transport}. The
parabolic equation \refeq{Transport} has the amplitude $\phi$ as the only unknown
function left, for which we intend to give an asymptotic expansion in the sequel. As
before we introduce a change of variables which adjusts the different time scales for the variables
of the tangent bundle, in particular:
\begin{equation}\label{eq:3.18}
  \psi(t,\xi,y) = t^\alpha \phi(t, t\xi+\tilde{x}(t),y+\tilde{y}(t))
\end{equation}
where the constant $\alpha$ is determined by consistency requirements. In terms of the
transformed function the transport equation takes the form
\begin{eqnarray}\label{eq:3.19}
\lefteqn{0=\frac{\partial \psi}{\partial t} - \frac{\alpha}t\psi
   + \frac 12 ~tr\left(g(t\xi + \tilde{x}) \frac{\partial^2 \Sigma}{\partial y^2}\right)
   + (g(t\xi+\tilde{x})\frac{\partial \Sigma}{\partial y},\frac{\partial \psi}{\partial y})
   -(g(\tilde{x})\tilde{q}, \frac{\partial \psi}{\partial y})}\\
&&-t^{-1}(\xi + G(t\xi + \tilde{x})(y+\tilde{y})
   - G(\tilde{x})\tilde{y},\frac{\partial \psi}{\partial \xi})
   + \frac 12 (\frac{\partial G}{\partial x}(t\xi+\tilde{x})(y+\tilde{y}),
             (y+\tilde{y})\frac{\partial \psi}{\partial y}) \nonumber
\end{eqnarray}
\begin{proposition}\label{prop:3.4}
There exists a unique $\alpha > 0$ such that there exists an asymptotic expansion for the
solution of \eqref{eq:3.3}
\begin{displaymath}
  \psi = \psi_0 + t\psi_1 + t^2\psi_2+\ldots
\end{displaymath}
where each $\psi_k,\ k\ge 0$, constitutes a polynomial in $(\xi,y)$ with degree strictly
smaller than k and $\psi_0$ is some constant. Moreover, the solution is unique up to a
constant multiplier.
\end{proposition}
\begin{proof}
Assuming the solution of \eqref{eq:3.19} to be of the form of an asymptotic expansion
(regular power series) the equation \eqref{eq:3.19} decomposes into separate equations for each order
in t, e.g. for the lowest order $(t^{-1})$ we find
\begin{displaymath}
  -\alpha\psi_0
  + \frac 12\psi_0 ~tr(g(x_0)\frac{\partial^2 \Sigma_{-1}}{\partial y^2}(x_0,y_0))
  = 0\ ,
\end{displaymath}
where we also have inserted the asymptotic expansion given by Proposition \ref{prop:3.1}.
This determines the parameter $\alpha$ in \eqref{eq:3.18} to be
\begin{equation}\label{eq:3.20}
  \alpha = \frac 12~tr(g(x_0)\frac{\partial^2 \Sigma_{-1}}{\partial y^2}(x_0,y_0))\quad .
\end{equation}
Recalling that we are working in normal coordinates and reminding the representation given for
$\Sigma_{-1}$ in Lemma \ref{lemma:3.2} we find $\alpha = 2n$. Comparing the coefficients of the
term $t^0$ reveals the following equation
\begin{displaymath}
  \psi_1 - \left(\xi+y, \frac{\partial \psi_1}{\partial \xi} \right)
         + (6\xi+4y, \frac{\partial \psi_1}{\partial y}) = 0
\end{displaymath}
the solution of which vanishes due to Lemma \ref{lemma:3.1}. In particular, this equation is
the homogeneous counterpart to \refeq{uNULL}. For the coefficients of the term $t^1$ we
get the following nonhomogeneous counterpart of the previous equation:
\begin{displaymath}
  \psi_2 - \left(\xi+y, \frac{\partial \psi_2}{\partial \xi}\right)
    + \left(6\xi+4y, \frac{\partial \psi_2}{\partial y}\right)
    +~tr\left(\frac 12 \frac{\partial^2 \Sigma_{-1}}{\partial y^2}
    + g^{kl}(\xi_k-y_k^0)(\xi_l-y_l^0)\right) = 0 .
\end{displaymath}
Decomposing the trace into its homogeneous polynomials of degree 0,1,2, respectively, we
may determine the solution (a polynomial of degree 2 in ($\xi,y,y^0$)) as indicated by
Lemma \ref{lemma:3.1}.
\begin{equation}
\begin{array}{l}
\psi =  1+{t}^{2}\left(-\left(-{\frac {3}{56}}\,g_{k,i}^{i,j}\xi_{{k}}-{
\frac {3}{56}}\,g_{k,i}^{j,i}\xi_{{k}}-{\frac {3}{56}}\,g_{k,j}^{i,i}\xi_{{k}} \right)\xi_{{j}}
+\frac{1}{6}\,\left({\frac {7}{10}}\,g_{i,j}^{i,k}y_{0,j}-\frac{1}{5}\,g_{k,j}^{i,i}y_{0,j}
\right.\right.\\[2mm]
\left.\left.%
\phantom{.}-{\frac {3}{56}}\,g_{k,i}^{i,j}y_{j}-{\frac {3}{56}}\,g_{k,i}^{j,i}y_{j}-
{\frac {3}{56}}\,g_{k,j}^{i,i}y_{{j}}-g_{i,i}^{k,l}y_{0,l}\right )y_{{k}}
+\frac{5}{6}\,\left ({\frac {7}{10}}\,g_{i,j}^{i,k}y_{0,j}-\frac{1}{5}\,g_{k,j}^{i,i}y_{0,j}
-{\frac {3}{56}}\,g_{k,i}^{i,j}y_{{j}}
\right.\right.\\[2mm]
\left. \left.%
\phantom{.}-{\frac {3}{56}}\,g_{k,i}^{j,i}y_{{j}}-{\frac {3}{56}}\,g_{k,j}^{i,i}y_{{j}}-g_{i,i}^{k,l}y_{0,l}\right )\xi_{{k}}
-\frac{1}{6}\,g_{i,i}^{k,l}y_{0,k}y_{{l}}-\frac{5}{6}\,g_{i,i}^{k,l}y_{0,k}\xi_{{l}}
+{\frac {1}{210}}\,g_{i,j}^{i,k}y_{{k}}y_{{j}}
\right.\\[2mm]
\phantom{.}+{\frac {3}{140}}\,g_{i,j}^{i,k}\xi_{{j}}y_{{k}}+{\frac {4}{35}}\,g_{i,j}^{i,k}\xi_{{j}}\xi_{{k}}+\frac{1}{30}\,g_{i,i}^{k,l}y_{{l}}y_{{k}}
+{\frac {3}{20}}\,g_{i,i}^{k,l}\xi_{{k}}y_{{l}}+\frac{4}{5}\,g_{i,i}^{k,l}\xi_{{k}}\xi_{{l}}
\\[2mm]
\left.
\phantom{.}-\frac{1}{30}\,g_{i,i}^{k,l}y_{0,k}y_{0,l}+\frac{1}{15}\,g_{i,j}^{i,k}y_{0,j}y_{0,k}-\frac{1}{30}\,g_{k,j}^{i,i}y_{0,k}y_{0,j}\right )+O(t^3)\ .  \hss
\end{array}
\label{phi:eqn}
\end{equation}
\end{proof}
The inverse transformation
$\phi(t,x,y) = \psi(t,\frac{x-\tilde{x}(t)}t, y-\tilde{y}(t))$ has
been done explicitly for an approximation of order 1 using MAPLE.
Applying this back transformation to $\psi$ (which is given by an asymptotic expansion in $t$) will
give back an asymptotic expansion in $t$. In particular, for $x_0=0$ we have
\begin{displaymath}
  \phi(t,x,y,y_0) = (1+ O(|x| + |y-y^0| + t)) t^{-\alpha}\quad .
%\mbox{n\"{a}chste Ordnung in $y-y_0$}
\end{displaymath}

Inserting the expansions for $S(t,x,y,x_0,y_0)$ and $\phi$ into the WKB Ansatz
\eqref{eq:3.1} reveals directly that the leading term is Gaussian.

\vskip 0.2cm
We now proceed to show that the error term can be expressed by a factor.
For an arbitrary m-dimensional system of ordinary differential equations
$\frac d{dt}{\cal{X}}(t,\alpha) = f(\cal{X},\alpha)$ depending on an m-dimensional
parameter $\alpha$ a result of Liouville states that the square matrix
${\cal J} = [\frac {\partial{\calX}} {\partial \alpha}]_{ij}$ is non degenerate on
some time interval and satisfies
$\frac d{dt}{\calJ} = {\calJ}~tr[\frac{\partial f}{\partial {\calX}}]$.
%This can be seen by induction.
%Beweis (Induktionsanfang) siehe Nebenrechnung
For a proof see e.g. \cite{FeM81}.
It is possible to determine the amplitude $\phi$ explicitly by applying this result to the system
\begin{displaymath}
  \frac d{dt}x = \frac{\partial H}{\partial p}(x,y,p,q) \quad\mbox{ and } \quad
  \frac d{dt}y = \frac{\partial H}{\partial q}(x,y,p,q)
\end{displaymath}
where $H$ is given in \eqref{Hamilton} and by using the equivalence of solutions
for the Hamilton and the Hamilton-Jacobi equations given by Propositions
\ref{prop:2.5} and \ref{cor:2.2ii}. We find
\begin{displaymath}
  p = \frac{\partial S}{\partial x}(t,x,y,x_0,y^0)\quad\mbox{ and } \quad
  q = \frac{\partial S}{\partial y}(t,x,y,x_0,y^0)
\end{displaymath}
to be dependent variables of $x,y$.
Hence along the characteristics $(X,Y)(t,x_0,y^0,q^0,p_0)$ of the Hamilton equation
corresponding to the Hamiltonian \eqref{Hamilton} the determinant
$J = det \frac{\partial (X,Y)}{\partial (q^0,p_0)}$ satisfies the equation
\begin{displaymath}
  \dot{J} = J~tr
  \left( \left[ \begin{array}{ll}
  \frac{\partial^2 H}{\partial p \partial x} & \frac{\partial^2 H}{\partial p \partial y}\\
  \frac{\partial^2 H}{\partial q \partial x} & \frac{\partial^2 H}{\partial q \partial y}\\
  \end{array} \right]
  +
  \left[ \begin{array}{ll}
  \frac{\partial^2 H}{\partial p^2} & \frac{\partial^2 H}{\partial p \partial q}\\
  \frac{\partial^2 H}{\partial q \partial p} & \frac{\partial^2 H}{\partial q^2}\\
  \end{array} \right]
  \left[ \begin{array}{ll}
  \frac{\partial^2 S}{\partial x^2} & \frac{\partial^2 S}{\partial x \partial y}\\
  \frac{\partial^2 S}{\partial y \partial x} & \frac{\partial^2 S}{\partial y^2}\\
  \end{array} \right] \right)
\end{displaymath}
where the second term in this sum corresponds to the dependent variables $p,\ q$.\\
By inserting the definition of the Hamiltonian \eqref{Hamilton} most partial derivatives of $H$
vanish leaving us with
\begin{equation}\label{eq:3.24a}
  \dot J = J  tr \left(-\frac{\partial G}{\partial x} y +  \frac{\partial G}{\partial x} y
                +  \frac{\partial^2 H}{\partial q^2}\frac{\partial^2 S}{\partial y^2}\right)\\
         = J tr \left(g(x) \frac{\partial^2 S}{\partial y^2}\right)
\end{equation}
which implies
\begin{equation}\label{eq:3.24b}
  \frac d{dt}J^{-\frac 12}
  = -\frac 12 J^{-\frac 12}~tr\left(g(x)\frac{\partial^2 S}{\partial y^2}\right)\ .
\end{equation}
On the other hand inserting the Hamilton equations into the transport equation \eqref{Transport}
reveals
\begin{displaymath}
  \frac d{dt} \phi
  + \frac 12 \phi~tr\left(\frac{\partial^2 H}{\partial q^2}
             \left(x,y,\frac{\partial S}{\partial x}, \frac{\partial S}{\partial y}\right)
             \frac{\partial^2 S}{\partial^2 (x,y)} \right) = 0
\end{displaymath}
or more explicitly
\begin{equation}\label{eq:3.25}
  \frac d{dt} \phi + \frac 12 \phi~tr\left(g(x)\frac{\partial^2 S}{\partial y^2}\right) = 0\quad .
\end{equation}
Comparing \eqref{eq:3.24b} and \eqref{eq:3.25} we find that $J^{-\frac 12}$ and $\phi$ are
identical up to a constant factor $\alpha$.

We now deal with the inhomogeneity in the PDE \refeq{uNULL}. Aiming at a multiplicative error
for an exponential approximation of WKB type respectively a  multiplicative error for
$u^0$ we need uniform estimates for the inhomogeneity in an appropriate neighbourhood of the
starting point $(x_0, y^0)$. This will be achieved by exchanging $\phi$ with $J{^-\frac 12}$
and using the estimates given in Corollary \ref{cor:2.2i} and \refeq{eq:2.1b}.
\begin{theorem}\label{thm:3.1} Let us assume that the assumptions of
Theorem~\ref{thm:2.3} are satisfied in the whole cotangent bundle.
Let an approximating solution of the form \eqref{eq:3.3} to \eqref{Greensfunktion}
satisfy the Hamilton Jacobi- \eqref{Hamilton-Jacobi} and transport equation
\eqref{Transport}. Then the inhomogenity $F_0$ \eqref{uNULL} is a multiplicative
error of the type
\begin{displaymath}
  \left[ h\frac{\partial}{\partial t} -
  H\left(x,y,-h\frac{\partial}{\partial x}, -h\frac{\partial}{\partial y}\right) \right] u^0
  = O(h^2) t^{2} u^0
\end{displaymath}
\end{theorem}
\begin{proof} Recalling the
definition of $F_0$ and the fact that $\phi$ and $J^{-\frac 12}$ are identical up
to a constant factor c.f. \eqref{eq:3.24b} and \eqref{eq:3.25} it suffices to
examine
\begin{displaymath}
  \frac{\partial^\nu}{(\partial x_i)^\nu} J^{-\frac 12}(t,x,y,x_0,y^0) \quad\mbox{and}\quad
  \frac{\partial^\nu} {(\partial y_i)^\nu}J^{-\frac 12}(t,x,y,x_0,y^0)
\end{displaymath}
for $\nu = 1,~2$ and $1\leq i \leq n$. Let $p_0 = p_0(t,x,y,x_0,y^0)$ and
$q_0 = q_0(t,x,y,x_0,y^0)$ satisfy the Hamilton Jacobi equation
\eqref{Hamilton-Jacobi}. Then we have along the characteristics
$(X,Y)(t,x_0,y^0,q^0,p_0)$ of the Hamiltonian given in \eqref{Hamilton} for
$\gamma = (\gamma_0,\gamma_1)$ with $\gamma_0 := x$, $\gamma_1 := y$
both in $\RR^d$ and $\zeta = (q^0,p_0)\in \RR^{2d}$
\begin{eqnarray}\label{eq:3.25a}
   \frac{\partial}{\partial \gamma_i} J^{-\frac 12}
&=& -\frac 12 J^{-\frac 12} \left(J^{-1}\frac{\partial J}{\partial \gamma_i}\right)
 = -\frac 12 J^{-\frac 12}~tr\left[ \left( \frac{\partial (X,Y)}{\partial \zeta} \right)^{-1}
      \frac{\partial}{\partial \gamma_i}\frac{\partial (X,Y)}{\partial \zeta}\right]\nonumber\\
&=&-\frac 12 J^{-\frac 12}\left(\frac{\partial^2 (X,Y)}{\partial\zeta\partial\zeta}\right)^k_{lm}
  \left( \frac{\partial (X,Y) }{\partial \gamma}\right)_{li}^{-1}
  \left( \frac{\partial (X,Y) }{\partial \zeta}\right)_{mk}^{-1}
\end{eqnarray}
since $\frac{\partial (X,Y) }{\partial \gamma_i} =
\frac{\partial (X,Y) }{\partial \zeta} \frac{\partial \zeta}{\partial \gamma_i}$ .
Analogously, by iteration
\begin{displaymath}
\frac{\partial^2}{\partial \gamma_i^2}J^{-\frac 1 2}
= - \frac 12 \frac{\partial}{\partial \gamma_i}\left(J^{-\frac 32}
                \frac{\partial}{\partial \gamma_i}J\right)
    = -\frac 12\left[ -\frac 32J^{-\frac 52}(\frac{\partial}{\partial \gamma_i}J)^2
          +  J^{-\frac 32}\frac{\partial^2}{\partial \gamma_i^2}J \right]
    = -\frac 12 J^{-\frac 12} K
\end{displaymath}
where
\begin{eqnarray*}
K:&=& -\frac 32 J^{-2}\frac{\partial}{\partial \gamma_i}J)^2
      + J^{-1}\frac{\partial^2 J }{\partial\zeta^l\partial\zeta^n}
               \frac{\partial \zeta^l \partial \zeta^n}{\partial \gamma_i^2}
      + J^{-1}\frac{\partial}{\partial \zeta^l}J
               \frac{\partial^2 \zeta^l}{\partial \gamma_i^2} \\
  &=& -\frac 32 \left(\frac{\partial (X,Y)}{\partial \zeta}\right)_{mk}^{-2}
       \left( \left(\frac{\partial^2 (X,Y)}{\partial \zeta^2}\right)_{lm}^k
       \left(\frac{\partial (X,Y)}{\partial \gamma}\right)_{li}^{-1}\right)^2  \\
  & &+~\left(\frac{\partial (X,Y)}{\partial \zeta}\right)_{mk}^{-1}
       \left(\frac{\partial^3 (X,Y)}{\partial \zeta^l \partial \zeta^m \partial \zeta^n}\right)
       \left(\frac{\partial (X,Y)}{\partial \gamma}\right)_{li}^{-1}
       \left(\frac{\partial (X,Y)}{\partial \zeta}\right)_{ni}^{-1} \\
  & &+~\left(\frac{\partial (X,Y)}{\partial \zeta}\right)_{km}^{-1}
       \left(\frac{\partial^2 (X,Y)}{\partial \zeta^l \partial \zeta^m}\right)^k
       \left(\frac{\partial^2 \zeta}{\partial \gamma_i^2}\right)^{lm}\ .\\
\end{eqnarray*}
Here the last factor may be rewritten as follows
$
  \frac{\partial^2 \zeta}{\partial \gamma_i^2}
  = -\left(\frac{\partial \gamma}{\partial \zeta}\right)_{\ell i}^{-2}
  \frac{\partial^2 \gamma_i}{\partial\zeta^l\partial\zeta^n}
  \left(\frac{\partial \gamma}{\partial \zeta}\right)_{in}^{-1}.
$
In order to simplify notations we replace $I_m$ by $m$, etc. in the next formula
where $I_m$ was defined in (\ref{eq:2.3}).
Applying Corrollary \ref{cor:2.2i} we find the following asymptotic behaviour for
$K$:
\begin{eqnarray*}
  \lefteqn{O(t^{2})+ O(t^{-(3-m-k)})O(t^{3\cdot 3-k-\ell-m-n})O(t^{-(3-\ell-i)})O(t^{-(3-n-i)})
  +O(t^{-(3-k-m)})} \\
  &\times&O(t^{3\cdot 3-k-\ell-m-n}t^{-(3-n-i)}t^{-(3-\ell -i)})
  -O(t^{-2(3-\ell -i)} t^{2\cdot 3-i-n-\ell}t^{2\cdot 3-m-k-\ell}t^{-(3-i-n)})                  \\
  &=& O(t^2i) + O(t^2i) + O(t^2i) \ .
\end{eqnarray*}
Recalling that $i=I_i=1$ for $y_i$, $1\le i\le d$, and inserting into \refeq{eq:3.25a}
reveals for the leading terms
\begin{equation} \label{eq:3.26}
   \frac{\partial }{\partial y_i} J^{-\frac 12}= J^{-\frac 12}O(t^1)
\end{equation}
and by iteration
\begin{equation} \label{eq:3.27}
  \frac{\partial^2 }{\partial y_i^2} J^{-\frac 12}= J^{-\frac 12}O(t^2)
\end{equation}
which finishes the proof of Theorem~\ref{thm:3.1}.
\end{proof}
Next we relax the assumption that Theorem \ref{thm:3.1}
hold in the whole cotangent  bundle by introducing a cutoff function in the
definition of the function $u_0$.
Let $\chi_D$ be a smooth function with values in $[0,1]$ such that $\chi_D$ vanishes
outside of the domain $D$ of the cotangent bundle, is 1 for all points $(x,y)$ in the
domain $D$ except in a neighbourhood of the boundary $\partial D$ of $D$, moreover, all
partial derivatives of $\chi_D$ until order one in $x$ and until order two in $y$ are
bounded for $t\in(0,t_0)$, $t_0 > 0$. Then the function
\begin{equation}\label{u-local}
  u^D = C(h) \chi_D(x-x_0,y-y_0)\phi(t,x,y,x_0,y_0)\exp[-S(t,x,y,x_0,y_0)/h]
\end{equation}
is globally well-defined for $t\leq t_0$ and satisfies
\begin{eqnarray}\label{eq:3.28}
h\frac{\partial u^D}{\partial t}\hspace{-2mm}&-&\hspace{-3mm}\frac{h^2}2~tr\left(g(x)
\frac{\partial^2 u^D}{\partial y_i \partial y_j}\right)
- h\left(G(x)y,\frac{\partial u^D}{\partial x}\right)
+ \frac h2 \left(\frac{\partial}{\partial x}(G(x)y,y),
\frac{\partial u^D}{\partial y}\right)\nonumber\\
  &=& -h^2\, C(h) \,F(t,x,y,x_0,y^0)
\end{eqnarray}
where the inhomogeneity is now given by
\begin{eqnarray*}
 F &=& tr\left( G(x)\frac{\partial^2 \phi}{\partial y_i \partial y_j}\right) \chi_D
               \exp\left[-\frac{S}{h}\right]
       + \frac{h^{-1}}{2}\left(tr\left(g(x)\frac{\partial \chi_D}{\partial y_j}
                                      \frac{\partial S}{\partial y_i}\right)\phi\right. \\
   &&\left.+~tr\left(g(x)\frac{\partial \chi_D}{\partial y_j}
                         \frac{\partial \phi}{\partial y_i}   \right)
           +~tr\left(g(x)\frac{\partial^2 \chi_D}{\partial y_j \partial y_i}\right)\phi\right)
                                                          \exp\left[-\frac{S}{h}\right]\\
   && +~h^{-2}\left(\left(G(x)y,\frac{\partial \chi_D}{\partial x}\right)
      -\frac{1}{2}\left(\frac{\partial}{\partial x}(G(x)y,y)
       \frac{\partial \chi_D}{\partial y}\right)\right)\phi\exp\left[-\frac{S}{h}\right]\ .
\end{eqnarray*}
In the sequel we are going to study the order of magnitude in $t$ and $h$ of the
inhomogeneity $F$. First we collect some general assumptions and facts. The
appropriate neighbourhood is the
polydisk $ D= B_r(x_0)\times B_{\frac rt}(y_0)$ in particular its closure $\bar D$.
The functions $\chi_D$, $\frac{\partial}{\partial x}\chi_D$, $\frac{\partial}{\partial y}\chi_D$
and $\frac{\partial^2}{\partial y^2}\chi_D$ are bounded on $D$, moreover, $g$ and
$G=g^{-1}$ are bounded since the manifold is compact. In normal coordinates
(actually $x_0=0$) $g(x)=id + G_{ij}x_ix_j+\ldots$ and
$G(x)=id - G_{ij}(x-x_0)_i(x-x_0)_j+\ldots$
which implies that $O(\frac{\partial}{\partial x}g) = O(x-x_0)$ = $O(-y_0 t)$ = $c\land r$
and  $O(\frac{\partial}{\partial x}G) = O(-(x-x_0))$ = $O(y_0 t)$ = $c\land r$.
Using the expansion given for $S$ in \refeq{eq:2.7} we find
\begin{eqnarray*}
  \frac{\partial S}{\partial y_i}
  &=& \frac 6{t^2}(1-\delta)(x-x_0)_i + \frac 2t(1-\delta)(y-y_0)_i\\
  &=& \frac 6{t^2}(1-\delta) c\land r + \frac 2t(1-\delta)(y_0 + t q_0 + y_0 + O(t))_i .
\end{eqnarray*}
this implies $O(\frac{\partial S}{\partial y_i})=O(t^{-2}) + O(t^{-2}) = O(t^{-2})$.
Using Proposition \ref{prop:2.5} we get similarly
$O(\frac{\partial S}{\partial y_i})=O(q_0)= O(t^{-2})$. Applying \refeq{eq:3.25a}
to $\frac{\partial \phi}{\partial y_i}$ and collecting terms Theorem \ref{thm:3.1}
is modified by the cutoff $\chi_D$ showing the following order of magnitude in $t$ and $h$
for the inhomogeneity $F$:
\begin{eqnarray*}
 O(t^2)u^D
   \!\! &+&\!\! \left(\frac {h^{-1}}2 ( O(t^{-2}) + O(t) + 1 )
   + h^{-2} (O(t^{-1}) + O(t^1) O(t^{-2})\! \right)\!\chi_{\Delta D_\varepsilon}\phi \exp\! \left[\! -\frac Sh\right] \\
\!\! &\le& \!\!O(t^2)u^D + \frac {h^{-1}}2 (O(t^{-2}) +1 + O(t)) \chi_{\Delta D_\varepsilon}
        O\left(e^{-\frac{\Omega_0-\delta}{ht{^-3}}}\right)
\quad.
\end{eqnarray*}
where $\chi_{\Delta D_\varepsilon}$ stands for $\frac{\partial}{\partial y}\chi_D$
as well as $\frac{\partial^2}{\partial y_i y_j}\chi_D$. The support
$\Delta D_\varepsilon$ of these functions is a
neighbourhood of the boundary $\partial D$ of width $\varepsilon$. In the above
estimate we used that for given $\delta>0$, $h>0$ one can choose $t_0$ such that
for $t\le t_0$ and $h\le h_0$
\begin{displaymath}
 12^n(2\pi t)^{-2n} \le \exp\left[\frac{\delta}{ht^3}\right],\qquad
 \frac 1{ht^3}\le \exp\left[ \frac{\delta}{ht^3}\right]\ .
\end{displaymath}
The constant $\Omega=\Omega_0 + \delta$ summarizes the constants
$\Omega_0=\min\{t^3 S(z,\xi)\vert
\norm{x-x_0}=r-\varepsilon\mbox{ and } \norm{y-y_0}=\frac{r-\varepsilon}{t}\}$ and
$\delta$. These calculations are summarized in the following proposition.

\begin{proposition}\label{prop:3.3}
Under the assumptions of Theorem \ref{thm:2.3} the inhomogeneity $F$ has the
asymptotic behaviour
$O(t^2) u^D + O(\exp [-\frac {\Omega }{ht^{-3}} ])\chi_{\Delta D_\varepsilon}$ for
$t,h>0$ and some positive constant $\Omega$.
\end{proposition}
After having studied the solution \refeq{u-local} of WKB shape we passed on to solving
the inhomogeneous equation \refeq{eq:3.28}. We now reconstruct the solution to the
original homogeneous equation using the well known DuHamel principle. Recalling the
classical results of Egorov, Hörmander, Oleinik and others we know that equation
\refeq{eq:1} is hypoelliptic since the generator $\cal L$ is of the form
\begin{displaymath}
 {\cal L} = \Sigma_{i=1}^n X_i^2 + X_0
\end{displaymath}
where $X_i= A_{ij}\frac{\ptl}{\ptl y_j}$, $1\le i,j\le n$ with $A^*A=G$,
$det A \ne 0$, and $X_0=(G(x)y,\frac {\ptl}{\ptl x})
- \frac 12 \left(\frac {\ptl}{\ptl x} (G(x)y,y),\frac {\ptl}{\ptl y}\right)$,
moreover, the Lie-brackets
$\left[X_i,(G(x)y,\frac {\ptl}{\ptl x})\right] = (A^*)^{-1}_{i\ell}\frac {\ptl}{\ptl x}$,
$1\le i\le n$, together with the vector fields $X_i$, $1\le i\le n$, form a basis
of the tangent space. In
particular it is shown cf. Hörmander \cite{Hor68}, Oleinik \cite{Ole71} that for any
initial distribution in $\cal D'$ and smooth inhomogeneity there exists a unique smooth
solution to \refeq{eq:1}. From this it immediately follows that the operators
$\left(e^{t\cal L}\right)_{t > 0}$ possess a smoothing kernel.
For the sake of a shorter notation we are using $z=(x,y)$ and $z_0=(x_0,y^0)$.
\begin{proposition}
For any domain $D$ in $\RR^{2d}$ and $C^2$ function $f$ with support in $D$ there exists
a unique smooth solution $v_t$ in $D$ satisfying the inhomogeneous equation
\begin{displaymath}
  \frac {\ptl}{\ptl t} v - {\cal L} v = f \ .
\end{displaymath}
In fact, $v$ is given by the integral formula
\begin{displaymath}
   v(t,z,z_0) = u(t,z,z_0)
                -\int_0^t\int_{\RR^{2n}} u(t-\tau,z,\xi)f(\tau,\xi,z_0)\, d\xi d\tau
\end{displaymath}
where $u$ satisfies the corresponding homogeneous equation (\ref{Greensfunktion}).
\end{proposition}
\begin{proof}
Uniqueness follows from the uniqueness of the corresponding homogeneous
equation \refeq{eq:1}. As for existence, applying the operator
$ \frac {\ptl}{\ptl t} - {\cal L}$ to the r.h.s. above yields
\begin{eqnarray*}
  \lefteqn{\left( \frac {\ptl}{\ptl t} - {\cal L}\right)v
 = \left( \frac {\ptl}{\ptl t} - {\cal L}\right) u(t,z,z_0)
  -  \left( \frac {\ptl}{\ptl t} - {\cal L}\right)
     \int_0^t\int_{\RR^{2n}} u(t-\tau,z,\xi)f(\tau,\xi,z_0)\, d\xi d\tau} \\
 &=& -\frac {\ptl}{\ptl t}
       \int_0^t\int_{\RR^{2n}} u(t-\tau,z,\xi)f(\tau,\xi,z_0)\, d\xi d\tau
     + {\cal L} \int_0^t\int_{\RR^{2n}} u(t-\tau,z,\xi)f(\tau,\xi,z_0)\, d\xi d\tau\\
 &=& \int_{\RR^{2n}} u(0,z,\xi)f(\tau,\xi,z_0)\, d\xi
  - \int_0^t\int_{\RR^{2n}} \frac {\ptl}{\ptl t} u(t-\tau,z,\xi)f(\tau,\xi,z_0)\, d\xi d\tau\\
  && + {\cal L}\int_0^t\int_{\RR^{2n}} u(t-\tau,z,\xi)f(t,\xi,z_0)\, d\xi d\tau\\
 &=&  f(t,z_0,z_0) \ .
\end{eqnarray*}
We have used that $u$ solves the homogeneous equation to eliminate terms.
Furthermore when differentiating the parameter integral w.r. to $t$ we have skipped the
regularization needed for the argument giving
$\int_{\RR^{2n}} \delta_z(\xi)f(t,\xi,z_0)\, d\xi$ and have exploited the cutoff
given by the function $f$ for obtaining the integrability of
$(\frac {\ptl}{\ptl t} u) f$.
\end{proof}
For the homogeneous equation \refeq{Greensfunktion} and its nonhomogeneous counterpart
\refeq{eq:3.28} the solutions are therefore related in the following way
\begin{displaymath}
  u^D = (1- h I_t)u\ ,
\end{displaymath}
with $u^D$ given by \refeq{u-local} and where the operator $I_t$ is defined by
$I_t u (t,z,z_0):= \int_0^t\int_{\RR^{2n}} u(t-\tau,z,\xi) F(\tau,\xi,z_0)\, d\xi d\tau$
with $F$ given in \refeq{eq:3.28}. In case the inverse operator exists, which is
equivalent to the convergence of its associated formal series $I+ \sum h^m I_t^m$, we
retrieve the solution to \refeq{Greensfunktion}.
\begin{theorem}\label{finalerror}
For small $t$ the fundamental solution $u$ to \refeq{Greensfunktion}
has the form
\begin{displaymath}
 u = u^D(1+O(ht)) + O(\exp [-\frac {\Omega }{ht^{-3}} ])
\end{displaymath}
with $u^D$ given by \refeq{u-local} and $\Omega$ as in Proposition \ref{prop:3.3}.
\end{theorem}

In order to meet the assumptions of Proposition \ref {prop:2.1}, Theorem \ref{thm:2.3},
and Proposition \ref{prop:3.3} we choose $D=B_r(x_0)\times B_{\frac rt}(y_0)$ for $r$
sufficiently small. Due to Proposition \ref{prop:3.3} the action functional $S$ satisfies
$\min_{\xi}\{S(t-\tau,z,\xi)+S(\tau,\xi,z_0)\}=S(t,z,z_0)$. Introducing a zero term, using that
$S$ is a convex function on $D$ with a unique minimum, and replacing
$S(t-\tau,z,\xi)+S(\tau,\xi,z_0)-\min_{\xi}\{S(t-\tau,z,\xi)+S(\tau,\xi,z_0)\}$
by its Taylor polynomial of order three at $\xi_0$ minimizing
$S(t-\tau,z,\xi)+S(\tau,\xi,z_0)$ we find:
\begin{eqnarray}\label{eq:3.29}
 \lefteqn{I_t u^D(z)
 \le C\int_0^t\int_{\RR^{2d}}  u^D(t-\tau,z,\xi)\left(\tau^2 u^D(\tau,\xi,z_0)
     + e^{-\frac \Omega{h\tau^3}}\chi_{\Delta D_\varepsilon}(\xi-z_0)\right)
         \, d\xi d\tau}\nonumber\\
 &\le&  e^{-\frac 1h S(t,z,z_0)}\!\!\int_0^t
      \!\!\frac{C(h)^2\tau^2}{(t-\tau)^{2d} \tau^{2d}}
     \!\!\int_{\RR^{2d}}\!\!
     e^{-\frac1h ((\Xi_{f,t-\tau}+\Xi_{b,\tau})(\xi-\xi_0),\xi-\xi_0)}
     \chi_D(\xi-z_0)\chi_D(z-\xi)\, d\xi\, d\tau\nonumber\\
 &&+\quad C(h) \!\!\int_0^t\int_{\RR^{2d}}\!\!(t-\tau)^{-2d}e^{-\frac 1h S(t-\tau,z,\xi)}\chi_D(z-\xi)
           e^{-\frac \Omega{h\tau^3}}\chi_{\Delta D\varepsilon}(\xi-z_0)
           \, d\xi\, d\tau\nonumber\\
 &\le& C C(h) t^{-2d} e^{-\frac 1h S(t,z,z_0)}t\chi_{D^{\bigstar 2}}(z-z_0)\nonumber \\
 && + C(h) \int_0^t e^{-\frac \Omega{h\tau^3}} e^{-M}
    \!\!\int_{\RR^{2d}}\!\!(t-\tau)^{-2d}\!\!
             e^{-\frac 1h (\Xi_{b,\tau}(z-\xi),z-\xi)}\chi_D(z-\xi)
             \chi_{D^{\bigstar 2}}(z-z_0)\, d\xi\, d\tau\nonumber\\
 &\le& C \left( C(h)^2 e^{-\frac 1h S(t,z,z_0)}t^{-2d+1}
       +  t e^{-\frac \Omega{ht^3}}\right) \chi_{D^{\bigstar 2}}(z-z_0)
 \end{eqnarray}
where $C(h)=\frac{12^d}{(4\pi ht)^{2d}}$ and
$\chi_{D^{\bigstar 2}}(z-z_0)$ is the characteristic function of
the poly disc
$D^{\bigstar 2}=\{(x,y)\in \RR^{2d}\vert\norm{x-x_0}\le 2r,\, \norm{y-y_0}\le \frac{2r}t \}$.
Moreover, above we have introduced the matrices:
$\Xi^{ij}_{b,t}=\frac{\partial^2 S(t,z,\xi)}{\partial y'_i\partial y'_j}
$, and
$\Xi^{ij}_{f,t}=\frac{\partial^2 S(t,\xi,z_0)}{\partial y'_i\partial y'_j}
$, $i\le i,j\le d$,
having non vanishing $2\times 2$-matrices on the diagonal
as given in \refeq{invcovariance}, where $\xi=(x',y')\in\RR^{2d}$.
The determinants of both matrices are given by $\frac{12^{d}}{\sigma^{2d}}$ with
$\sigma\in\{t-\tau,\tau\}$. Their sum equals $t^{2d}\tau^{-2d}(t-\tau)^{-2d}$.
All constants are included in the universal constant $C$.
Finally $\Omega$ is given in Proposition \ref{prop:3.3}. We have used
$\chi_D(\xi-z_0)\chi_D(z-\xi)\le \chi_{D^{\bigstar 2}}(z-z_0)$ and that
the function $\exp(-\frac{\Omega}{ht^3})$ is isotone having the image $[0,1)$.
The $m^{th}$ order term $I^m_t$ follows by iteration. For the sake of a clear
representation we introduce the following notations:
\begin{eqnarray*}
  f(t,z,z_0) &=& t^{-2d}\chi_D(z-z_0)C(h) exp(-S(t,z,z_0)/h), \\
  g_m(t,z,z_0) &=& \chi_{D^{\bigstar m}}(z-z_0)\exp\left[-\frac{\Omega}{ht^3}\right]
\end{eqnarray*}
with $z=(x,y)\in\RR^{2d}$, $z_0=(x_0,y_0)\in\RR^{2d}$ and
$\chi_{D^{\bigstar m}}(z-z_0)$ is the characteristic function of
the polydisc
$D^{\bigstar m}=\{(x,y)\in \RR^{2d}\vert
              \norm{x-x_0}\le mr,\, \norm{y-y_0}\le \frac{mr}t \}$.
In terms of the new notation
inequality \refeq{eq:3.29}
reads $f\ast f (t,z,z_0)\le C C(h)  f(t,z,z_0)$ together with
$f\ast C g_1(t,z,z_0)\le  C  tg_2(t,z,z_0)$ for some
constant $C$.
As to the higher order convolutions we apply:
\begin{eqnarray*}
 \lefteqn{\norm{f\ast g_m (t,z,z_0)}   =  \int_0^t \int_{\RR^{2d}}\!\!
     \chi_{D^{\bigstar m}}(z-\xi)  e^{-\frac{\Omega}{h(t-\tau)^3}}
     \frac{C(h)}{ \tau^{2d}}
      \norm{e^{-\frac 1h S(\tau,\xi,z_0)}}  \chi_D(\xi-z_0)   \, d\xi\, d\tau}\\
 &\le&  \chi_{D^{\bigstar m+1}}(z-z_0) \int_0^t e^{-\frac{\Omega}{h(t-\tau)^3}}
        e^{-M} \int_{\RR^{2d}}  \frac{C(h)}{ \tau^{2d}}
        e^{(\Xi_{f,t}(\xi-z_0),\xi-z_0)}   \chi_D(\xi-z_0) \, d\xi\, d\tau\\
 &\le& g_{m+1}(t,z,z_0) e^{-M} t.
\end{eqnarray*}
where we have used that $S$ is locally convex satisfying
$S(t,z,z_0)\ge M + (\Xi_{f,t}(\xi-z_0),\xi-z_0)$ on $D$. Moreover, there holds
\begin{eqnarray*}
 \lefteqn{ g_1\ast g_m (t,z,z_0) = \int_0^t \int_{\RR^{2d}}\!\!
        \chi_{D^{\bigstar m}}(z-\xi)  e^{-\frac{\Omega}{h(t-\tau)^3}}
        e^{-\frac{\Omega}{h\tau^3}} \chi_D(\xi-z_0)   \, d\xi\, d\tau }\\
  &\le& \chi_{D^{\bigstar m+1}}(z-z_0)
         e^{-\frac{\Omega}{ht^3}} \int_0^t   e^{-\frac{\Omega}{h\tau^3}}
        \int_{\RR^{2d}}  \chi_D(\xi-z_0) \, d\xi\, d\tau\\
  &\le& \pi^{n} \Gamma\left(\frac d2 +1\right)^2 r^{d}  \chi_{D^{\bigstar m+1}}(z-z_0)
         e^{-\frac{\Omega}{ht^3}} \int_0^t \left({\frac rt}\right)^{d}
         e^{-\frac{\Omega}{h\tau^3}} d\tau
  \le t b(r) g_{m+1}(t,z,z_0).
\end{eqnarray*}
where we have estimated
$t^{-n}e^{-\frac{\Omega}{h\tau^3}}
  \le t^{-n}e^{-\frac{\Omega}{ht^3}}
  \le e^{-\frac{\Omega-\delta}{ht^3}}\le 1
$
and $b(r)$ summarizes the constants.
Hence by induction for the $m^{th}$ order term we find
\begin{eqnarray*}
 h^{m-1}I_t^{\bigotimes m-1} u^D(z)
 &\le&  C C(h) (C h b(r)t)^{m-1}(f+2g_2+\ldots 2^{m-1} g_m)\\
 &\le& C C(h) (C h b(r)t)^{m-1} (f+g_m).
 \end{eqnarray*}
This is applied when studying convergence, i.e. estimating
$\norm{(\sum_{0}^\infty h^m I_t^{\bigotimes m})u^D-u}(t,x,y,x_0,y^0)$.
Since $C h b(r)t < 1$ for $h<h_0$ and $t<t_0$, $h_0,t_0$ sufficiently small, we
get
\begin{displaymath}
C(h)2C \sum_{\ell=m}^\infty (C h b(r)t)^{\ell} g_\ell
    \le C(h)\frac{(2C h b(r)t)^{m}}{1- (2C h_0 b(r)t_0)^{m}}  e^{-\frac{\Omega}{ht^3}}
    \le \frac{(2C h b(r)t)^{m}}{1- (2C h_0 b(r)t_0)^{m}} e^{-\frac{\Omega-\delta}{ht^3}}
\end{displaymath}
where the ratio tends to zero for $m\rightarrow\infty$. Moreover, we also get that
$C(h)f \sum_{\ell=m}^\infty (C h b(r)t)^{\ell}
\le C(h)f \frac{(2C h b(r)t)^{m}}{1- (2C h_0 b(r)t_0)^{m}}$ tends to zero for
$m\rightarrow\infty$. Hence the series $\sum_0^{\infty} h^{m}I_t^{\bigotimes m} $
converges uniformly.
Hence we have shown that the solution $u$ of the homogeneous system is given by a
convergent series on any compact set if $h$ is small.
\vskip 0.3cm
{\bf Remark}
The result is local in as much as $u^D$ is localized in $D$. But the $\Omega$
in the above estimates in Theorem \ref{finalerror} is uniform. Since the manifold is compact this result is
enough for getting global results. This concludes the studies of the heat kernel.
Our result is related to the construction
of a generalized heat kernel on a bundle over a compact manifold in \cite{BGV91}
concedrning the structure of the underlying space, however the diffusion is
degenerate.\\

\section{A trace formula}

In this last section we use the uniform bounds derived above to derive a trace
formula, i.e. an asymptotic expansion for the trace of the semigroup describing the
evolution of the solution to equation \refeq{Greensfunktion}.
Let us explain the most important facts. Let $H$ and $H'$ be
(separable, infinite dimensional) Hilbert spaces, and choose orthonormal bases
$(e_i)_{i\in\NN}$ and $(e_j)_{j\in\NN}$ in $H$ and $H'$, respectively. For any
bounded operator $A: H\rightarrow H'$ the quantity
$\Norm{A}_{HS}=\sum_{i,j} \langle A e_i,e_j\rangle$ is
independent of the choice of bases in $H$ and $H'$, c.f. \cite{Sim79}\cite{Roe98}.
An operator if it is finite such that $\Norm{A}_{HS}<\infty$ is called an {\em Hilbert Schmidt operator}
and $\Norm{A}_{HS}$ is called its {\em Hilbert Schmidt norm}. The Hilbert Schmidt norm
is induced by a scalar product
$\langle A,B\rangle=
               \sum_{i,j}{\bar{\langle A e_i,e_j\rangle}}\langle Be_j,e_i\rangle$.
A bounded operator is of {\em trace class} if there are
Hilbert Schmidt operators
$A$ and $B$ on $H$ with $T=AB$. Its {\em trace} $Tr(T)$ is defined to be the Hilbert
Schmidt inner product $\langle A^*,B\rangle_{HS}$. The trace $Tr(T)$ depends in
fact only on $T$ i.e. it is independent of the choice of $A$ and $B$.

In \cite{BGV91} chapter 2 a covariant approach to a generalized Laplacian with
corresponding heat kernel on bundles over a compact orientable Riemannian manifold $M$
is developed. In particular, the scalar product between the space of compactly
supported sections of the bundle $\mathcal{E}$ and the space of sections of
$\mathcal{E}^*\otimes\norm{\Lambda}$ is established, where $\norm{\Lambda}$ is the
density bundle. The result given in Theorem 12.60 in \cite{CFKS87} can be
adopted. Since the generator ${\cal L}$ given by \refeq{Greensfunktion}
considered in this paper is nonsymmetric we will,
however, mainly exploit the existence of a kernel instead of working with the
operators. Since the kernel is smoothing on the space of distributions according
to \cite{Hor68} and \cite{Ole71} and since it constitutes the transistion probability of
a diffusion the operators $\exp[t {\cal L}]$, $t>0$,
the corresponding Hilbert Schmidt norm is given by
\begin{eqnarray*}
 \Norm{e^{t {\cal L}}}_{HS} &=& \int_{T^*M} u(x,y,x,y)\,dxdy
   \le \int_{M}\int_{TM} u^D(x,y,x,y)(1+O(t))dy\,dx\\
   &\le& \mbox{vol}(M) \int_{TM} e^{\frac 1h 6y^2}(1+O(t))dy <\infty
\end{eqnarray*}
is finite, see e.g. c.f. \cite{Sim79}.
Exploiting the semigroup property of a transition density
we see that the operators of the semi group are in fact of trace class
\begin{eqnarray*}
  Tr \exp[t {\cal L}] &=& \int u(t,x,y,x,y)\, dx\,dy
  =\int\!\!\!\int u({\frac t2},x,y,x',y')u({\frac t2},x',y',x,y)\, dx'dy'dx dy\\
  &=&\langle\exp[\frac t2 {\cal L}],\exp[\frac t2 {\cal L}]\rangle.
\end{eqnarray*}

Inserting the expansion we derived in Section 4 we find:
\begin{theorem} Under the hypothesis of Proposition \ref{prop:2.1} we have the
following asymptotic expansion with respect to a small parameter $t$.
\begin{displaymath}
Tr \exp[t {\cal L}]= \int_{T^*M} u(t,x,y,x,y) dx\,dy  =  {12^{n\over 2}\over (2\pi)^n t^{2n}}
(\mbox{vol} M + R t^3+a_4t^4+\cdots+a_\nu t^\nu +O(t^{\nu+1}) ) ,
\end{displaymath}
where $\mbox{vol} M$ is the volume of $M$ and the first non trivial coefficient $R$ is
the scalar curvature. All further coefficients $a_j$, $j=4,\ldots$ are invariant
quantities of $M$ which can be computed by the method indicated before.
\end{theorem}
\begin{center}{\bf Acknowledgements}\end{center}
\noindent The authors would like to thank Zdislav Brze\'zniak and David Elworthy for
stimulating discussions. For his work with the manuskript we are indebted to
Andrei Khrennikov.
The second author is particularly grateful to Remi L\'eandre for introducing her
to the theory of small parameter estimates and gratefully acknowledges financial
support by ''Profilen matematisk modellering'' at Växjö University.
The authors are deeply indebted to Paul Fischer for his help with programming in
MAPLE which was essential for the computations.

\bibliography{ASTRID}
\end{document}